# Local Linear Convergence Analysis of Primal–Dual Splitting Methods

Jingwei Liang*, Jalal Fadili† and Gabriel Peyré‡

**Abstract.** In this paper, we study the local linear convergence properties of a versatile class of Primal–Dual splitting methods for minimizing composite non-smooth convex optimization problems. Under the assumption that the non-smooth components of the problem are partly smooth relative to smooth manifolds, we present a unified local convergence analysis framework for these Primal–Dual splitting methods. More precisely, in our framework we first show that (i) the sequences generated by Primal–Dual splitting methods identify a pair of primal and dual smooth manifolds in a finite number of iteration, and then (ii) enter a local linear convergence regime, which is for instance characterized in terms of the structure of the underlying active smooth manifolds. We also show how our results for Primal–Dual splitting specialize to cover existing one on Forward–Backward splitting and Douglas–Rachford splitting/ADMM (alternating direction methods of multipliers). Moreover, based on these obtained local convergence analysis result, several practical acceleration techniques for the class of Primal–Dual splitting methods are discussed. To exemplify the usefulness of the obtained result, we consider several concrete numerical experiments arising from applicative fields including signal/image processing, inverse problems and machine learning, etc. The demonstration not only verify the local linear convergence behaviour of Primal–Dual splitting methods, but also the insights on how to accelerate them in practice.

## 1 Introduction

### 1.1 Composed optimization problem

In various fields such as inverse problems, signal and image processing, statistics and machine learning *etc.*, many problems are (eventually) formulated as structured optimization problems (see Section 6 for some specific examples). A typical example of these optimization problems, given in its primal form, reads

$$\min_{x\in\mathbb{R}^n} R(x) + F(x) + (J \mathbin{\overset{+}{\vee}} G)(Lx), \tag{$\mathcal{P}_{\mathrm{P}}$}$$

where $(J \mathbin{\overset{+}{\vee}} G)(\cdot) \stackrel{\mathrm{def}}{=} \inf_{v\in\mathbb{R}^m} J(\cdot) + G(\cdot - v)$ denotes the infimal convolution of $J$ and $G$. Throughout, we assume the following:

- (**A.1**) $R, F \in \Gamma_0(\mathbb{R}^n)$ with $\Gamma_0(\mathbb{R}^n)$ being the class of proper convex and lower semi-continuous functions on $\mathbb{R}^n$, and $\nabla F$ is $(1/\beta_F)$-Lipschitz continuous for some $\beta_F > 0$,
- (**A.2**) $J, G \in \Gamma_0(\mathbb{R}^m)$, and $G$ is $\beta_G$-strongly convex for some $\beta_G > 0$.
- (**A.3**) $L : \mathbb{R}^n \to \mathbb{R}^m$ is a linear mapping.

*DAMTP, University of Cambridge, UK, E-mail: jl993@cam.ac.uk

†Normandie Université, ENSICAEN, CNRS, GREYC, France, E-mail: Jalal.Fadili@ensicaen.fr

‡CNRS, DMA, ENS Paris, France, E-mail: Gabriel.Peyre@ens.fr

**(A.4)** $0 \in \mathrm{ran}\big(\partial R + \nabla F + L^*(\partial J \Box \partial G)L\big)$, where $\partial J \Box \partial G \stackrel{\text{def}}{=} (\partial J^{-1} + \partial G^{-1})^{-1}$ is the parallel sum of the subdifferential operators $\partial J$ and $\partial G$, and $\mathrm{ran}(\cdot)$ denotes the range of a set-valued operator. See Remark 3.2 for the reasoning of this condition.

The main difficulties encountered to solve such a problem are that the objective function is non-smooth, the presence of the linear operator $L$ and the infimal convolution. Consider also the Fenchel-Rockafellar dual problem [43]

$$\min_{v \in \mathbb{R}^m} J^*(v) + G^*(v) + (R^* \,\overset{+}{\vee}\, F^*)(-L^* v). \qquad (\mathcal{P}_{\mathrm{D}})$$

The classical Kuhn-Tucker theory asserts that a pair $(x^\star, v^\star) \in \mathbb{R}^n \times \mathbb{R}^m$ solves $(\mathcal{P}_{\mathrm{P}})$-$(\mathcal{P}_{\mathrm{D}})$ if it satisfies the monotone inclusion

$$0 \in \begin{bmatrix} \partial R & L^* \\ -L & \partial J^* \end{bmatrix} \begin{pmatrix} x^\star \\ v^\star \end{pmatrix} + \begin{bmatrix} \nabla F & 0 \\ 0 & \nabla G^* \end{bmatrix} \begin{pmatrix} x^\star \\ v^\star \end{pmatrix}, \qquad (1.1)$$

One observes that in (1.1), the composition by the linear operator and the infimal convolution have been decoupled, hence opening the door to achieve full splitting. This is a key property that is used by all Primal–Dual algorithms that we are about to review. In turn, solving (1.1) provides a pair of points that are solutions to $(\mathcal{P}_{\mathrm{P}})$ and $(\mathcal{P}_{\mathrm{D}})$ respectively.

More complex forms of $(\mathcal{P}_{\mathrm{P}})$ involving for instance a sum of infimal convolutions can be tackled in a similar way using a product space trick, as we will see in Section 5.

### 1.2 Primal–Dual splitting methods

Primal–Dual splitting methods to solve more or less complex variants of $(\mathcal{P}_{\mathrm{P}})$-$(\mathcal{P}_{\mathrm{D}})$ have witnessed a recent wave of interest in the literature [13, 11, 50, 18, 29, 21, 16]. All these approaches achieve full splitting, they involve the resolvents of $R$ and $J^*$, the gradients of $F$ and $G^*$ and the linear operator $L$, all separately at various points in the course of iteration. For instance, building on the seminal work of [3], the now-popular scheme proposed in [13] solves $(\mathcal{P}_{\mathrm{P}})$-$(\mathcal{P}_{\mathrm{D}})$ with $F = G^* = 0$. The authors in [29] have shown that the Primal–Dual splitting method of [13] can be seen as a proximal point algorithm (PPA) in $\mathbb{R}^n \times \mathbb{R}^m$ endowed with a suitable norm. Exploiting the same idea, the author in [21] considered $(\mathcal{P}_{\mathrm{P}})$ with $G^* = 0$, and proposed an iterative scheme which can be interpreted as a Forward–Backward (FB) splitting again with an appropriately renormed space. This idea is further extended in [50] to solve more complex problems such as that in $(\mathcal{P}_{\mathrm{P}})$. A variable metric version was proposed in [19]. Motivated by the structure of (1.1), [11] and [18] proposed a Forward–Backward-Forward scheme [48] to solve it.

In this paper, we will focus the unrelaxed Primal–Dual splitting method summarized in Algorithm 1. This scheme covers that of [13, 50, 29, 21, 16]. Though we omit the details here for brevity, our analysis and conclusions carry through to the method proposed in [11, 18].

---

**Algorithm 1:** A Primal–Dual splitting method

**Initial**: Choose $\gamma_R, \gamma_J > 0$ and $\theta \in [-1, +\infty[$. For $k = 0$, $x_0 \in \mathbb{R}^n$, $v_0 \in \mathbb{R}^m$;

**repeat**

$$\begin{aligned} x_{k+1} &= \mathrm{prox}_{\gamma_R R}(x_k - \gamma_R \nabla F(x_k) - \gamma_R L^* v_k), \\ \bar{x}_{k+1} &= x_{k+1} + \theta(x_{k+1} - x_k), \\ v_{k+1} &= \mathrm{prox}_{\gamma_J J^*}(v_k - \gamma_J \nabla G^*(v_k) + \gamma_J L \bar{x}_{k+1}), \end{aligned} \qquad (1.2)$$

$k = k + 1$;

**until** *convergence*;

---

**Remark 1.1.**

(i) Algorithm 1 is somehow an interesting extension to the literature given the choice of $\theta$ that we advocate. Indeed, the range of $\theta$ is $[-1, +\infty[$ is larger than the one proposed in [29] (*i.e.* $\theta \in [-1, 1]$). It encompasses the Primal–Dual splitting method proposed in [50] when $\theta = 1$, and the one in [21] when moreover $G^* = 0$. When both $F = G^* = 0$, it reduces to the Primal–Dual splitting method proposed in [13, 29].

(ii) It can also be verified that Algorithm 1 covers the Forward–Backward (FB) splitting [39] ($J^* = G^* = 0$), Douglas–Rachford (DR) splitting [25] (if $F = G^* = 0$, $L = \mathrm{Id}$ and $\gamma_R = 1/\gamma_J$) as special cases; see Section 4 for a discussion or [13] and references therein for more details. Exploiting this relation, in Section 4, we build connections with the results provided in [36, 37] for FB-type methods and DR/ADMM. It also should be noted that, DR splitting is the limiting case of the Primal–Dual splitting [13], and the global convergence result of Primal–Dual splitting does not apply to DR.

### 1.3 Contributions

In the literature, most studies on the convergence rate of Primal–Dual splitting methods mainly focus on the global behaviour [13, 23, 35, 9, 14]. For instance, it is now known that the (partial) duality gap decreases sublinearly at the rate $O(1/k)$, which can be accelerated to $O(1/k^2)$ [13] under strong convexity of either the primal or the dual problem. The iterates converge globally linearly if both the primal and dual problems are strongly convex [13, 9], or locally linearly under certain regularity assumptions [35]. However, in practice, local linear convergence of the sequences generated by Algorithm 1 has been observed for many problems with the absence of strong convexity (as confirmed by our numerical experiments in Section 6). None of the existing theoretical analysis was able to explain this behaviour so far. Providing the theoretical underpinnings of this local behaviour is the main goal pursued in this paper. Our main findings can be summarized as follows.

**Finite time activity identification** For Algorithm 1, let $(x^\star, v^\star)$ be a Kuhn-Tucker pair, *i.e.* a solution of (1.1). Under a non-degeneracy condition, and provided that both $R$ and $J^*$ are partly smooth relative to $C^2$-smooth manifolds, respectively $\mathcal{M}^R_{x^\star}$ and $\mathcal{M}^{J^*}_{v^\star}$ near $x^\star$ and $v^\star$ (see Definition 2.10), we show that the generated primal-dual sequence $\{(x_k, v_k)\}_{k\in\mathbb{N}}$ which converges to $(x^\star, v^\star)$ will identify in finite time the manifold $\mathcal{M}^R_{x^\star} \times \mathcal{M}^{J^*}_{v^\star}$ (see Theorem 3.3). In plain words, this means that after a finite number of iterations, say $K$, we have $x_k \in \mathcal{M}^R_{x^\star}$ and $v_k \in \mathcal{M}^{J^*}_{v^\star}$ for all $k \geq K$.

**Local linear convergence** Capitalizing on this finite identification result, we first show in Proposition 3.6 that the globally non-linear iteration (1.2) locally linearizes along the identified smooth manifolds, then we deduce that the convergence of the sequence becomes locally linear (see Theorem 3.11). The rate of linear convergence is characterized precisely based on the properties of the identified partly smooth manifolds and the involved linear operator $L$.

Moreover, when $F = G^* = 0, L = \mathrm{Id}$ and $R, J^*$ are locally polyhedral around $(x^\star, v^\star)$, we show that the convergence rate is parameterised by the cosine of the largest principal angle (yet smaller than $\pi/2$, see Definition 2.6) between the tangent spaces of the two manifolds at $(x^\star, v^\star)$ (see Lemma 3.8). This builds a clear connection between the results in this paper and those we drew in our previous works on DR and ADMM [38, 37].

### 1.4 Related work

For the past few years, an increasing attention has been paid to investigate the local linear convergence of first-order proximal splitting methods in absence of strong convexity. This has been done for instance for FB-type splitting [10, 27, 2, 30, 47], and DR/ADMM [24, 8, 4] for special objective functions. In our previous work [34, 36, 38, 37], based on the powerful framework provided by partial smoothness, we unified all the above-mentioned works and provide even stronger claims.

To the best of our knowledge, we are aware of only one recent paper [46] which investigated finite identification and local linear convergence of a Primal–Dual splitting method to solve a very special instance of ($\mathcal{P}_{\mathrm{P}}$). More precisely, they assumed $R$ to be gauge, $F = \frac{1}{2}\|\cdot\|^2$ (hence strong convexity of the primal problem), $G^* = 0$ and $J$ the indicator function of a point. Our work goes much beyond this limited case. It also deepens our current understanding of local behaviour of proximal splitting algorithms by complementing the picture we started in [36, 37] for FB and DR splitting.

**Paper organization** The rest of the paper is organized as follows. Some useful prerequisites, including partial smoothness, are collected in Section 2. The main contributions of this paper, *i.e.* finite time activity identification and local linear convergence of Primal–Dual splitting under partial smoothness are the core of Section 3. Several discussions on the obtained result are delivered in Section 4. Section 5 extends the results to the case of more than one infimal convolution. In Section 6, we report various numerical experiments to support our theoretical findings. All the proofs of the main results are collected in Section A.

# 2 Preliminaries

Throughout the paper, $\mathbb{N}$ is the set of non-negative integers, $\mathbb{R}^n$ is a $n$-dimensional real Euclidean space equipped with scalar product $\langle\cdot,\,\cdot\rangle$ and norm $\|\cdot\|$. $\mathrm{Id}_n$ denotes the identity operator on $\mathbb{R}^n$, where $n$ will be dropped if the dimension is clear from the context.

**Sets** For a nonempty convex set $C \subset \mathbb{R}^n$, denote $\mathrm{aff}(C)$ its affine hull, and $\mathrm{par}(C)$ the smallest subspace parallel to $\mathrm{aff}(C)$. Denote $\iota_C$ the indicator function of $C$, and $\mathrm{P}_C$ the orthogonal projection operator onto the set.

**Functions** The sub-differential of a function $R \in \Gamma_0(\mathbb{R}^n)$ is a set-valued operator,

$$\partial R : \mathbb{R}^n \rightrightarrows \mathbb{R}^n,\ x \mapsto \big\{g \in \mathbb{R}^n | R(x') \geq R(x) + \langle g,\, x' - x\rangle, \forall x' \in \mathbb{R}^n\big\}, \tag{2.1}$$

which is maximal monotone (see Definition 2.2). For $R \in \Gamma_0(\mathbb{R}^n)$, the proximity operator of $\gamma R$ is

$$\mathrm{prox}_{\gamma R}(x) \stackrel{\mathrm{def}}{=} \mathrm{argmin}_{z\in\mathbb{R}^n} \frac{1}{2}\|z - x\|^2 + \gamma R(z). \tag{2.2}$$

Given a function $J \in \Gamma_0(\mathbb{R}^m)$, its *Legendre-Fenchel conjugate* is defined as

$$J^*(y) = \sup_{v\in\mathbb{R}^m} \langle v,\, y\rangle - J(v). \tag{2.3}$$

We have $J^* \in \Gamma_0(\mathbb{R}^m)$ and $J^{**} = J$.

**Lemma 2.1 (Moreau Identity [41]).** *Let $J \in \Gamma_0(\mathbb{R}^n)$, then for any $x \in \mathbb{R}^n$ and $\gamma \in ]0, +\infty[$,*

$$x = \mathrm{prox}_{\gamma J}(x) + \gamma \mathrm{prox}_{J^*/\gamma}(x/\gamma). \tag{2.4}$$

Using the Moreau identity, it is straightforward to see that the update of $v_k$ in Algorithm 1 can be deduced also from $\mathrm{prox}_{J/\gamma_J}$.

**Operators** Given a set-valued mapping $A : \mathbb{R}^n \rightrightarrows \mathbb{R}^n$, its range is $\mathrm{ran}(A) = \{y \in \mathbb{R}^n : \exists x \in \mathbb{R}^n \text{ s.t. } y \in A(x)\}$, and graph is $\mathrm{gph}(A) \stackrel{\text{def}}{=} \{(x, u) \in \mathbb{R}^n \times \mathbb{R}^n | u \in A(x)\}$.

**Definition 2.2 (Monotone operator).** A set-valued mapping $A : \mathbb{R}^n \rightrightarrows \mathbb{R}^n$ is said to be monotone if,

$$\langle x_1 - x_2,\, v_1 - v_2 \rangle \geq 0,\ \forall (x_1, v_1) \in \mathrm{gph}(A) \ \text{ and } \ (x_2, v_2) \in \mathrm{gph}(A). \tag{2.5}$$

It is moreover maximal monotone if $\mathrm{gph}(A)$ can not be contained in the graph of any other monotone operator. If there exists an $\kappa > 0$ such that

$$\langle x_1 - x_2,\, v_1 - v_2 \rangle \geq \kappa \|x_1 - x_2\|^2,$$

then $A$ is called *strongly monotone*, and its *inverse* $A^{-1}$ is $\kappa$-cocoercive.

Let $\beta \in ]0, +\infty[$, $B : \mathbb{R}^n \to \mathbb{R}^n$, then $B$ is $\beta$-cocoercive if the following holds

$$\langle B(x_1) - B(x_2),\, x_1 - x_2 \rangle \geq \beta \|B(x_1) - B(x_2)\|^2,\ \forall x_1, x_2 \in \mathbb{R}^n, \tag{2.6}$$

which implies that $B$ is $\beta^{-1}$-Lipschitz continuous.

For a maximal monotone operator $A$, $(\mathrm{Id} + A)^{-1}$ is its resolvent. It is known that for $R \in \Gamma_0(\mathbb{R}^n)$ and $\gamma > 0$, $\mathrm{prox}_{\gamma R} = (\mathrm{Id} + \gamma \partial R)^{-1}$ [6, Example 23.3].

**Definition 2.3 (Non-expansive operator).** An operator $\mathcal{F} : \mathbb{R}^n \to \mathbb{R}^n$ is non-expansive if

$$\forall x, y \in \mathbb{R}^n,\ \|\mathcal{F}(x) - \mathcal{F}(y)\| \leq \|x - y\|.$$

For any $\alpha \in ]0, 1[$, $\mathcal{F}$ is called $\alpha$-averaged if there exists a non-expansive operator $\mathcal{F}'$ such that $\mathcal{F} = \alpha \mathcal{F}' + (1 - \alpha)\mathrm{Id}$.

The class of $\alpha$-averaged operators is closed under relaxation, convex combination and composition [6, 20]. In particular when $\alpha = \frac{1}{2}$, $\mathcal{F}$ is called *firmly non-expansive*, and several properties of it are collected in the following lemma.

**Lemma 2.4.** *Let $\mathcal{F} : \mathbb{R}^n \to \mathbb{R}^n$. Then the following statements are equivalent:*

(i) *$\mathcal{F}$ is firmly non-expansive;*
(ii) *$\mathrm{Id} - \mathcal{F}$ is firmly non-expansive;*
(iii) *$2\mathcal{F} - \mathrm{Id}$ is non-expansive;*
(iv) *Given any $\lambda \in ]0, 2]$, $(1 - \lambda)\mathrm{Id} + \lambda \mathcal{F}$ is $\frac{\lambda}{2}$-averaged;*
(v) *$\mathcal{F}$ is the resolvent of a maximal monotone operator $A : \mathbb{R}^n \rightrightarrows \mathbb{R}^n$.*

**Proof.** (i)⇔(ii)⇔(iii) follow from [6, Proposition 4.2, Corollary 4.29], (i)⇔(iv) is [6, Corollary 4.29], and (i)⇔(v) is [6, Corollary 23.8]. □

Let $\mathcal{S}(\mathbb{R}^n) = \{\mathcal{V} \in \mathbb{R}^{n \times n} | \mathcal{V}^T = \mathcal{V}\}$ the set of symmetric matrices acting on $\mathbb{R}^n$. The Loewner partial ordering on $\mathcal{S}(\mathbb{R}^n)$ is defined as

$$\forall \mathcal{V}_1, \mathcal{V}_2 \in \mathcal{S}(\mathbb{R}^n),\ \ \mathcal{V}_1 \succcurlyeq \mathcal{V}_2 \iff \forall x \in \mathbb{R}^n, \langle \mathcal{V}_1 x,\, x \rangle \geq \langle \mathcal{V}_2 x,\, x \rangle.$$

Given any positive constant $\nu \in ]0, +\infty[$, define $\mathcal{S}_\nu$ as

$$\mathcal{S}_\nu \stackrel{\text{def}}{=} \big\{\mathcal{V} \in \mathcal{S}(\mathbb{R}^n) : \mathcal{V} \succcurlyeq \nu \mathrm{Id}\big\}, \tag{2.7}$$

*i.e.* the set of *symmetric positive definite matrices* whose eigenvalues are bounded below by $\nu$. For any $\mathcal{V} \in \mathcal{S}_\nu$, define the following induced scalar product and norm,

$$\langle x,\, x' \rangle_{\mathcal{V}} = \langle x,\, \mathcal{V} x' \rangle,\ \ \|x\|_{\mathcal{V}} = \sqrt{\langle x,\, \mathcal{V} x \rangle},\ \forall x, x' \in \mathbb{R}^n.$$

By endowing the Euclidean space $\mathbb{R}^n$ with the above scalar product and norm, we obtain the Hilbert space which is denoted by $\mathbb{R}^n_{\mathcal{V}}$.

**Lemma 2.5.** *Let the operators $A : \mathbb{R}^n \rightrightarrows \mathbb{R}^n$ be maximal monotone, $B : \mathbb{R}^n \to \mathbb{R}^n$ be $\beta$-cocoercive, and $\mathcal{V} \in \mathcal{S}_\nu$. Then for $\gamma \in ]0, 2\beta\nu[$,*

(i) $(\mathrm{Id} + \gamma\mathcal{V}^{-1}A)^{-1} : \mathbb{R}^n_{\mathcal{V}} \to \mathbb{R}^n_{\mathcal{V}}$ *is* firmly non-expansive*;*

(ii) $(\mathrm{Id} - \gamma\mathcal{V}^{-1}B) : \mathbb{R}^n_{\mathcal{V}} \to \mathbb{R}^n_{\mathcal{V}}$ *is* $\frac{\gamma}{2\beta\nu}$-averaged non-expansive*;*

(iii) *The operator* $(\mathrm{Id} + \gamma\mathcal{V}^{-1}A)^{-1}(\mathrm{Id} - \gamma\mathcal{V}^{-1}B) : \mathbb{R}^n_{\mathcal{V}} \to \mathbb{R}^n_{\mathcal{V}}$ *is* $\frac{2\beta\nu}{4\beta\nu-\gamma}$-averaged non-expansive.

**Proof.**

(i) See [19, Lemma 3.7];

(ii) Since $B : \mathbb{R}^n \to \mathbb{R}^n$ is $\beta$-cocoercive, given any $x, x' \in \mathbb{R}^n$, we have

$$\begin{aligned}\langle x - x',\, \mathcal{V}^{-1}B(x) - \mathcal{V}^{-1}B(x')\rangle_{\mathcal{V}} &\geq \beta\|B(x) - B(x')\|^2 \\ &= \beta\mathcal{V}\langle \mathcal{V}^{-1}B(x) - \mathcal{V}^{-1}B(x'),\, \mathcal{V}^{-1}B(x) - \mathcal{V}^{-1}B(x')\rangle_{\mathcal{V}} \\ &\geq \beta\nu\|\mathcal{V}^{-1}B(x) - \mathcal{V}^{-1}B(x')\|^2_{\mathcal{V}},\end{aligned}$$

which means $\mathcal{V}^{-1}B : \mathbb{R}^n_{\mathcal{V}} \to \mathbb{R}^n_{\mathcal{V}}$ is $(\beta\nu)$-cocoercive. The rest of the proof follows [6, Proposition 4.33].

(iii) See [42, Theorem 3]. □

## 2.1 Angles between subspaces

In this part we introduce the notions of principal angles and Friedrichs angle between two subspaces $T_1$ and $T_2$ of $\mathbb{R}^n$. Without the loss of generality, assume $1 \leq p \stackrel{\text{def}}{=} \dim(T_1) \leq q \stackrel{\text{def}}{=} \dim(T_2) \leq n - 1$.

**Definition 2.6 (Principal angles).** The principal angles $\theta_k \in [0, \frac{\pi}{2}]$, $k = 1, ..., p$ between subspaces $T_1$ and $T_2$ are defined by, with $u_0 = v_0 \stackrel{\text{def}}{=} 0$, and

$$\begin{aligned}\cos(\theta_k) \stackrel{\text{def}}{=} \langle u_k,\, v_k\rangle = \max\langle u,\, v\rangle \ \text{ s.t. } \ &u \in T_1, v \in T_2, \|u\| = 1, \|v\| = 1, \\ &\langle u,\, u_i\rangle = \langle v,\, v_i\rangle = 0,\ i = 0, \cdots, k-1.\end{aligned}$$

The principal angles $\theta_k$ are unique with $0 \leq \theta_1 \leq \theta_2 \leq \cdots \leq \theta_p \leq \pi/2$.

**Definition 2.7 (Friedrichs angle).** The Friedrichs angle $\theta_F \in ]0, \frac{\pi}{2}]$ between $T_1$ and $T_2$ is

$$\cos\big(\theta_F(T_1, T_2)\big) \stackrel{\text{def}}{=} \max\langle u,\, v\rangle \,\text{s.t.}\, u \in T_1 \cap (T_1 \cap T_2)^\perp, \|u\| = 1,\ v \in T_2 \cap (T_1 \cap T_2)^\perp, \|v\| = 1.$$

The following lemma shows the relation between the Friedrichs and principal angles whose proof can be found in [5, Proposition 3.3].

**Lemma 2.8 (Principal angles and Friedrichs angle).** *The Friedrichs angle is exactly $\theta_{d+1}$ where $d \stackrel{\text{def}}{=} \dim(T_1 \cap T_2)$. Moreover, $\theta_F(T_1, T_2) > 0$.*

**Remark 2.9.** The principal angles can be obtained by the singular value decomposition (SVD). For instance, let $X \in \mathbb{R}^{n\times p}$ and $Y \in \mathbb{R}^{n\times q}$ form the orthonormal basis for the subspaces $T_1$ and $T_2$ respectively. Let $U\Sigma V^T$ be the SVD of $X^TY \in \mathbb{R}^{p\times q}$, then $\cos(\theta_k) = \sigma_k,\ k = 1, 2, \ldots, p$ and $\sigma_k$ corresponds to the $k$'th largest singular value in $\Sigma$.

## 2.2 Partial smoothness

In this section, we introduce the notion of partial smoothness, which lays the foundation of our local convergence rate analysis. The concept of partial smoothness is first introduced in [33]. This concept, as well as that of identifiable surfaces [51], captures the essential features of the geometry of non-smoothness which are along the so-called active/identifiable manifold. For convex functions, a closely related idea is developed in [32]. Loosely speaking, a partly smooth function behaves *smoothly* as we move along the identifiable submanifold, and *sharply* if we move transversal to the manifold. In fact, the behaviour of the function and of its minimizers depend essentially on its restriction to this manifold, hence offering a powerful framework for algorithmic and sensitivity analysis theory.

Let $\mathcal{M}$ be a $C^2$-smooth embedded submanifold of $\mathbb{R}^n$ around a point $x$. To lighten the notation, henceforth we shall state $C^2$-manifold instead of $C^2$-smooth embedded submanifold of $\mathbb{R}^n$. The natural embedding of a submanifold $\mathcal{M}$ into $\mathbb{R}^n$ permits to define a Riemannian structure on $\mathcal{M}$, and we simply say $\mathcal{M}$ is a Riemannian manifold. $\mathcal{T}_{\mathcal{M}}(x)$ denotes the tangent space to $\mathcal{M}$ at any point near $x$ in $\mathcal{M}$. More materials on manifolds are given in Section A.1.

Below we present the definition of partly smooth function associated to functions in $\Gamma_0(\mathbb{R}^n)$.

**Definition 2.10 (Partly smooth function).** Let $R \in \Gamma_0(\mathbb{R}^n)$, and $x \in \mathbb{R}^n$ such that $\partial R(x) \neq \emptyset$. $R$ is then said to be *partly smooth* at $x$ relative to a set $\mathcal{M}$ containing $x$ if

(i) **Smoothness**: $\mathcal{M}$ is a $C^2$-manifold around $x$, $R$ restricted to $\mathcal{M}$ is $C^2$ around $x$;

(ii) **Sharpness**: The tangent space $\mathcal{T}_{\mathcal{M}}(x)$ coincides with $T_x \stackrel{\text{def}}{=} \mathrm{par}(\partial R(x))^{\perp}$;

(iii) **Continuity**: The set-valued mapping $\partial R$ is continuous at $x$ relative to $\mathcal{M}$.

The class of partly smooth functions at $x$ relative to $\mathcal{M}$ is denoted as $\mathrm{PSF}_x(\mathcal{M})$.

Capitalizing on the results of [33], it can be shown that under transversality assumptions, the set of partly smooth functions is closed under addition and pre-composition by a linear operator. Moreover, absolutely permutation-invariant convex and partly smooth functions of the singular values of a real matrix, *i.e.* spectral functions, are convex and partly smooth spectral functions of the matrix [22]. Popular examples of partly smooth functions are summarized in Section 6 whose details can be found in [36].

The next lemma gives expressions of the Riemannian gradient and Hessian (see Section A.1 for definitions) of a partly smooth function.

**Lemma 2.11.** *If $R \in \mathrm{PSF}_x(\mathcal{M})$, then for any $x' \in \mathcal{M}$ near $x$,*

$$\nabla_{\mathcal{M}} R(x') = \mathrm{P}_{T_{x'}}(\partial R(x')).$$

*In turn, for all $h \in T_{x'}$,*

$$\nabla^2_{\mathcal{M}} R(x')h = \mathrm{P}_{T_{x'}} \nabla^2 \widetilde{R}(x')h + \mathfrak{W}_{x'}\big(h, \mathrm{P}_{T_{x'}^{\perp}} \nabla \widetilde{R}(x')\big),$$

*where $\widetilde{R}$ is any smooth extension (representative) of $R$ on $\mathcal{M}$, and $\mathfrak{W}_x(\cdot,\cdot) : T_x \times T_x^{\perp} \to T_x$ is the Weingarten map of $\mathcal{M}$ at $x$.*

**Proof.** See [36, Fact 3.3]. □

# 3 Local linear convergence of Primal–Dual splitting methods

In this section, we present the main result of the paper, the local linear convergence analysis of Primal–Dual splitting methods. We first present the finite activity identification of the sequence $(x_k, v_k)$ generated by the methods, from which we further show that the fixed-point iteration of Primal–Dual splitting methods locally can be linearized, and the linear convergence follows naturally.

### 3.1 Finite activity identification

Let us first recall the result from [50, 17], that under a proper renorming, Algorithm 1 can be written as FB splitting. Let $\theta = 1$, from the definition of the proximity operator (2.2), we have that (1.2) is equivalent to the following monotone inclusion

$$- \begin{bmatrix} \nabla F & 0 \\ 0 & \nabla G^* \end{bmatrix} \begin{pmatrix} x_k \\ v_k \end{pmatrix} \in \begin{bmatrix} \partial R & L^* \\ -L & \partial J^* \end{bmatrix} \begin{pmatrix} x_{k+1} \\ v_{k+1} \end{pmatrix} + \begin{bmatrix} \mathrm{Id}_n/\gamma_R & -L^* \\ -L & \mathrm{Id}_m/\gamma_J \end{bmatrix} \begin{pmatrix} x_{k+1} - x_k \\ v_{k+1} - v_k \end{pmatrix}. \tag{3.1}$$

Define the product space $\mathcal{K} = \mathbb{R}^n \times \mathbb{R}^m$, and let $\mathbf{Id}$ be the identity operator on $\mathcal{K}$. Define the following variable and operators

$$\boldsymbol{z}_k \stackrel{\text{def}}{=} \begin{pmatrix} x_k \\ v_k \end{pmatrix}, \ \boldsymbol{A} \stackrel{\text{def}}{=} \begin{bmatrix} \partial R & L^* \\ -L & \partial J^* \end{bmatrix}, \ \boldsymbol{B} \stackrel{\text{def}}{=} \begin{bmatrix} \nabla F & 0 \\ 0 & \nabla G^* \end{bmatrix}, \ \boldsymbol{\mathcal{V}} \stackrel{\text{def}}{=} \begin{bmatrix} \mathrm{Id}_n/\gamma_R & -L^* \\ -L & \mathrm{Id}_m/\gamma_J \end{bmatrix}. \tag{3.2}$$

It is easy to verify that $\boldsymbol{A}$ is maximal monotone [11], $\boldsymbol{B}$ is $\min\{\beta_F, \beta_G\}$-cocoercive. For $\boldsymbol{\mathcal{V}}$, denote $\nu = (1 - \sqrt{\gamma_J \gamma_R \|L\|^2}) \min\{\frac{1}{\gamma_J}, \frac{1}{\gamma_R}\}$, then $\boldsymbol{\mathcal{V}}$ is symmetric and $\nu$-positive definite [50, 19]. Define $\boldsymbol{\mathcal{K}}_{\boldsymbol{\mathcal{V}}}$ the Hilbert space induced by $\boldsymbol{\mathcal{V}}$.

Now (3.1) can be reformulated as

$$\boldsymbol{z}_{k+1} = (\boldsymbol{\mathcal{V}} + \boldsymbol{A})^{-1}(\boldsymbol{\mathcal{V}} - \boldsymbol{B})\boldsymbol{z}_k = (\mathbf{Id} + \boldsymbol{\mathcal{V}}^{-1}\boldsymbol{A})^{-1}(\mathbf{Id} - \boldsymbol{\mathcal{V}}^{-1}\boldsymbol{B})\boldsymbol{z}_k. \tag{3.3}$$

Clearly, (3.3) is the FB splitting method on $\boldsymbol{\mathcal{K}}_{\boldsymbol{\mathcal{V}}}$ [50]. When $F = G^* = 0$, it reduces to the metric Proximal Point Algorithm discussed in [13, 29].

Before presenting the finite time activity identification under partial smoothness, we first recall the global convergence of Algorithm 1.

**Lemma 3.1 (Convergence of Algorithm 1).** *Consider Algorithm 1 under assumptions* (**A.1**)-(**A.4**). *Let $\theta = 1$ and choose $\gamma_R, \gamma_J$ such that*

$$2 \min\{\beta_F, \beta_G\} \min\Big\{\frac{1}{\gamma_J}, \frac{1}{\gamma_R}\Big\}\big(1 - \sqrt{\gamma_J \gamma_R \|L\|^2}\big) > 1, \tag{3.4}$$

*then there exists a Kuhn-Tucker pair $(x^\star, v^\star)$ such that $x^\star$ solves ($\mathcal{P}_{\mathrm{P}}$), $v^\star$ solves ($\mathcal{P}_{\mathrm{D}}$), and $(x_k, v_k) \to (x^\star, v^\star)$.*

**Proof.** See [50, Corollary 4.2]. □

**Remark 3.2.**

(i) Assumption (**A.4**) is important to ensure existence of Kuhn-Tucker pairs. There are sufficient conditions which ensure that (**A.4**) can be satisfied. For instance, assuming that ($\mathcal{P}_{\mathrm{P}}$) has at least one solution and some classical domain qualification condition is satisfied (see *e.g.* [18, Proposition 4.3]), assumption (**A.4**) can be shown to be in force.

(ii) It is obvious from (3.4) that $\gamma_J \gamma_R \|L\|^2 < 1$, which is also the condition needed in [13] for convergence in the special case $F = G^* = 0$. The convergence condition in [21] differs from (3.4), however, $\gamma_J \gamma_R \|L\|^2 < 1$ still is a key condition. The values of $\gamma_J, \gamma_R$ can also be made varying along iterations, and convergence of the iteration remains under the rule provided in [19]. However, for the sake of brevity, we omit the details of this case here.

(iii) Lemma 3.1 addresses global convergence of the iterates provided by Algorithm 1 only for the case $\theta = 1$. For the choices $\theta \in [-1, 1[\cup]1, +\infty[$, so far the corresponding convergence of the iteration cannot be obtained directly, and a correction step as proposed in [29] for $\theta \in [-1, 1[$ is needed so that the iteration is a contraction. Unfortunately, such a correction step leads to a new iterative scheme, not simply (1.2) itself, see [29] for more details. In a very recent paper [12], the authors also proved the convergence of the Primal–Dual splitting method of [50] for the case of $\theta \in [-1, 1]$ with a proper modification of the iterates. Since the main focus of this work is to investigate local convergence behaviour, the analysis of global convergence of Algorithm 1 for any $\theta \in [-1, +\infty[$ is beyond the scope of this paper. Thus, we will only consider the case $\theta = 1$ in our analysis. Nevertheless, as we will see later, locally $\theta > 1$ gives faster convergence rate compared to the choice $\theta \in [-1, 1]$. This points out a future direction of research to design new Primal–Dual splitting methods.

**Theorem 3.3 (Finite activity identification).** *Consider Algorithm 1 under assumptions (**A.1**)-(**A.4**). Let $\theta = 1$ and choose $\gamma_R, \gamma_J$ according to Lemma 3.1. Thus $(x_k, v_k) \to (x^\star, v^\star)$, where $(x^\star, v^\star)$ is a Kuhn-Tucker pair that solves $(\mathcal{P}_{\mathrm{P}})$-$(\mathcal{P}_{\mathrm{D}})$. If moreover $R \in \mathrm{PSF}_{x^\star}(\mathcal{M}^R_{x^\star})$ and $J^* \in \mathrm{PSF}_{v^\star}(\mathcal{M}^{J^*}_{v^\star})$, and the non-degeneracy condition holds*

$$
\begin{aligned}
-L^*v^\star - \nabla F(x^\star) &\in \mathrm{ri}\big(\partial R(x^\star)\big), \\
Lx^\star - \nabla G^*(v^\star) &\in \mathrm{ri}\big(\partial J^*(v^\star)\big).
\end{aligned}
\tag{ND}
$$

*Then,*

(i) *$\exists K \in \mathbb{N}$ such that for all $k \geq K$,*

$$
(x_k, v_k) \in \mathcal{M}^R_{x^\star} \times \mathcal{M}^{J^*}_{v^\star}.
$$

(ii) *Moreover,*

(a) *if $\mathcal{M}^R_{x^\star} = x^\star + T^R_{x^\star}$, then $T^R_{x_k} = T^R_{x^\star}$ and $\bar{x}_k \in \mathcal{M}^R_{x^\star}$ hold for $k > K$.*

(b) *If $\mathcal{M}^{J^*}_{v^\star} = v^\star + T^{J^*}_{v^\star}$, then $T^{J^*}_{v_k} = T^{J^*}_{v^\star}$ holds for $k > K$.*

(c) *If $R$ is locally polyhedral around $x^\star$, then $\forall k \geq K$, $x_k \in \mathcal{M}^R_{x^\star} = x^\star + T^R_{x^\star}$, $T^R_{x_k} = T^R_{x^\star}$, $\nabla_{\mathcal{M}^R_{x^\star}} R(x_k) = \nabla_{\mathcal{M}^R_{x^\star}} R(x^\star)$, and $\nabla^2_{\mathcal{M}^R_{x^\star}} R(x_k) = 0$.*

(d) *If $J^*$ is locally polyhedral around $v^\star$, then $\forall k \geq K$, $v_k \in \mathcal{M}^{J^*}_{v^\star} = v^\star + T^{J^*}_{v^\star}$, $T^{J^*}_{v_k} = T^{J^*}_{v^\star}$, $\nabla_{\mathcal{M}^{J^*}_{v^\star}} J^*(v_k) = \nabla_{\mathcal{M}^{J^*}_{v^\star}} J^*(v^\star)$, and $\nabla^2_{\mathcal{M}^{J^*}_{v^\star}} J^*(v_k) = 0$.*

See Section A.2 for the proof.

**Remark 3.4.**

(i) The non-degeneracy condition (ND) is a strengthened version of (1.1).

(ii) In general, we have no identification guarantees for $x_k$ and $v_k$ if the proximity operators are computed with errors, even if they are summable, in which case one can still prove global convergence. The reason behind this is that in the exact case, under condition (ND), the proximal mapping of the partly smooth function $R$ and that of its restriction to $\mathcal{M}^R_{x^\star}$ locally agree nearby $x^\star$ (and similarly for $J^*$ and $v^\star$). This property can be easily violated if approximate proximal mappings are involved, see [36] for an example.

(iii) Theorem 3.3 only states the existence of $K$ after which the identification of the sequences happens, and no bounds are available. In [36, 37], lower bounds of $K$ for the FB and DR splitting methods are provided, and similar lower bounds can be obtained here for the Primal–Dual splitting methods. Since such lower-bounds are only of theoretical interest, we decided to skip the corresponding details here and refer the reader to [36, 37].

### 3.2 Locally linearized iteration

Relying on the identification result, now we are able to show that the globally nonlinear fixed-point iteration (3.3) can be locally linearized along the manifold $\mathcal{M}^R_{x^\star} \times \mathcal{M}^{J^*}_{v^\star}$. As a result, the convergence rate of the iteration essentially boils down to analysing the spectral properties of the matrix obtained in the linearized iteration.

Given a Kuhn-Tucker pair $(x^\star, v^\star)$, define the following two functions

$$\overline{R}(x) \stackrel{\text{def}}{=} R(x) + \langle x,\, L^* v^\star + \nabla F(x^\star)\rangle, \;\; \overline{J^*}(y) \stackrel{\text{def}}{=} J^*(y) - \langle y,\, Lx^\star - \nabla G^*(v^\star)\rangle. \tag{3.5}$$

We have the following lemma.

**Lemma 3.5.** *Let $(x^\star, v^\star)$ be a Kuhn-Tucker pair such that $R \in \mathrm{PSF}_{x^\star}(\mathcal{M}^R_{x^\star})$, $J^* \in \mathrm{PSF}_{x^\star}(\mathcal{M}^{J^*}_{v^\star})$. Denote the Riemannian Hessians of $\overline{R}$ and $\overline{J^*}$ as*

$$H_{\overline{R}} \stackrel{\text{def}}{=} \gamma_R \mathrm{P}_{T^R_{x^\star}} \nabla^2_{\mathcal{M}^R_{x^\star}} \overline{R}(x^\star) \mathrm{P}_{T^R_{x^\star}} \;\; \text{and} \;\; H_{\overline{J^*}} \stackrel{\text{def}}{=} \gamma_J \mathrm{P}_{T^{J^*}_{v^\star}} \nabla^2_{\mathcal{M}^{J^*}_{v^\star}} \overline{J^*}(v^\star) \mathrm{P}_{T^{J^*}_{v^\star}}. \tag{3.6}$$

*Then $H_{\overline{R}}$ and $H_{\overline{J^*}}$ are symmetric positive semi-definite under either of the following conditions:*

(i) (ND) *holds.*

(ii) *$\mathcal{M}^R_{x^\star}$ and $\mathcal{M}^{J^*}_{v^\star}$ are affine subspaces.*

*Define,*

$$W_{\overline{R}} \stackrel{\text{def}}{=} (\mathrm{Id}_n + H_{\overline{R}})^{-1} \;\; \text{and} \;\; W_{\overline{J^*}} \stackrel{\text{def}}{=} (\mathrm{Id}_m + H_{\overline{J^*}})^{-1}, \tag{3.7}$$

*then both $W_{\overline{R}}$ and $W_{\overline{J^*}}$ are firmly non-expansive.*

**Proof.** See [37, Lemma 6.1]. □

In addition to (A.1) and (A.2), in the rest of the paper, we assume that $F$ and $G^*$ locally are $C^2$-smooth around $x^\star$ and $v^\star$ respectively. Now define the restricted Hessians of $F$ and $G^*$,

$$H_F \stackrel{\text{def}}{=} \mathrm{P}_{T^R_{x^\star}} \nabla^2 F(x^\star) \mathrm{P}_{T^R_{x^\star}} \;\; \text{and} \;\; H_{G^*} \stackrel{\text{def}}{=} \mathrm{P}_{T^{J^*}_{v^\star}} \nabla^2 G^*(v^\star) \mathrm{P}_{T^{J^*}_{v^\star}}. \tag{3.8}$$

Denote $\overline{H}_F \stackrel{\text{def}}{=} \mathrm{Id}_n - \gamma_R H_F, \overline{H}_{G^*} \stackrel{\text{def}}{=} \mathrm{Id}_m - \gamma_J H_{G^*}, \overline{L} \stackrel{\text{def}}{=} \mathrm{P}_{T^{J^*}_{v^\star}} L \mathrm{P}_{T^R_{x^\star}}$ and

$$M_{\mathrm{PD}} \stackrel{\text{def}}{=} \begin{bmatrix} W_{\overline{R}} \overline{H}_F & -\gamma_R W_{\overline{R}} \overline{L}^* \\ \gamma_J (1+\theta) W_{\overline{J^*}} \overline{L} W_{\overline{R}} \overline{H}_F - \theta \gamma_J W_{\overline{J^*}} \overline{L} & W_{\overline{J^*}} \overline{H}_{G^*} - \gamma_R \gamma_J (1+\theta) W_{\overline{J^*}} \overline{L} W_{\overline{R}} \overline{L}^* \end{bmatrix}. \tag{3.9}$$

We have the following proposition.

**Proposition 3.6 (Local linearized iteration).** *Suppose that Algorithm 1 is run under the conditions of Theorem 3.3. Then for all $k$ large enough, we have*

$$\boldsymbol{z}_{k+1} - \boldsymbol{z}^\star = M_{\mathrm{PD}} (\boldsymbol{z}_k - \boldsymbol{z}^\star) + o(\|\boldsymbol{z}_k - \boldsymbol{z}^\star\|). \tag{3.10}$$

See Section A.2 for the proof.

**Remark 3.7.**

(i) For the case of iteration varying $(\gamma_J, \gamma_R)$, yielding $\{(\gamma_{J,k}, \gamma_{R,k})\}_k$ according to the result of [36], (3.10) remains true if these parameters are converging to some constants such that condition (3.4) still holds.

(ii) Taking $\overline{H}_{G^*} = \mathrm{Id}_m$ (*i.e.* $G^* = 0$) in (3.9), one gets the linearized iteration associated to the Primal–Dual splitting method of [21]. If we further let $\overline{H}_F = \mathrm{Id}_n$, this will correspond to the linearized version of the method in [13].

Now we need to study the spectral properties of $M_{\mathrm{PD}}$. Let $p \stackrel{\text{def}}{=} \dim(T^R_{x^\star}), q \stackrel{\text{def}}{=} \dim(T^{J^*}_{v^\star})$ be the dimensions of the tangent spaces $T^R_{x^\star}$ and $T^{J^*}_{v^\star}$ respectively, define $S^R_{x^\star} = (T^R_{x^\star})^\perp$ and $S^{J^*}_{v^\star} = (T^{J^*}_{v^\star})^\perp$. Assume that $q \geq p$ (alternative situations are discussed in Remark 3.9). Let

$$\overline{L} = X \Sigma_{\overline{L}} Y^*$$

be the singular value decomposition of $\overline{L}$, and denote the rank of $\overline{L}$ as $l \stackrel{\text{def}}{=} \mathrm{rank}(\overline{L})$, clearly, we have $l \leq p$. Denote $M^k_{\mathrm{PD}}$ the $k$-th power of $M_{\mathrm{PD}}$.

**Lemma 3.8 (Convergence property of $M_{\text{PD}}$).** *The following holds for the matrix $M_{\text{PD}}$ in (3.9).*

(i) *If $\theta = 1$, then $M_{\text{PD}}$ is convergent to $M_{\text{PD}}^\infty$, i.e.*

$$M_{\text{PD}}^\infty \overset{\text{def}}{=} \lim_{k\to\infty} M_{\text{PD}}^k, \tag{3.11}$$

*which is the projection operator onto the set of fixed-points of $M_{\text{PD}}$. Moreover,*

$$\forall k \in \mathbb{N},\ M_{\text{PD}}^k - M_{\text{PD}}^\infty = (M_{\text{PD}} - M_{\text{PD}}^\infty)^k \ \text{ and }\ \rho(M_{\text{PD}} - M_{\text{PD}}^\infty) < 1.$$

*Given any $\rho \in ]\rho(M_{\text{PD}} - M_{\text{PD}}^\infty), 1[$, there is $K$ large enough such that for all $k \geq K$,*

$$\|M_{\text{PD}}^k - M_{\text{PD}}^\infty\| = O(\rho^k).$$

(ii) *If $F = G^* = 0$, and $R, J^*$ are locally polyhedral around $(x^\star, v^\star)$. Then given any $\theta \in ]0, 1]$, $M_{\text{PD}}$ is convergent with*

$$M_{\text{PD}}^\infty = \begin{bmatrix} Y & \\ & X \end{bmatrix} \begin{bmatrix} 0_l & & & \\ & \text{Id}_{n-l} & & \\ & & 0_l & \\ & & & \text{Id}_{m-l} \end{bmatrix} \begin{bmatrix} Y^* & \\ & X^* \end{bmatrix}. \tag{3.12}$$

*Moreover, all the eigenvalues of $M_{\text{PD}} - M_{\text{PD}}^\infty$ are complex with the spectral radius*

$$\rho(M_{\text{PD}} - M_{\text{PD}}^\infty) = \sqrt{1 - \theta\gamma_R\gamma_J\sigma_{\min}^2} < 1, \tag{3.13}$$

*where $\sigma_{\min}$ is the smallest non-zero singular value of $\overline{L}$.*

**Remark 3.9.** Here we discuss in short other possible cases of (3.12) when $F = G^* = 0$ and $R, J^*$ are locally polyhedral around $(x^\star, v^\star)$.

(i) When $L = \text{Id}$, then $\overline{L} = \mathrm{P}_{T_{v^\star}^{J^*}} \mathrm{P}_{T_{x^\star}^R}$ and $\sigma_{\min}$ stands for the cosine value of the biggest principal angle (yet strictly smaller than $\pi/2$) between tangent spaces $T_{x^\star}^R$ and $T_{v^\star}^{J^*}$.

(ii) For the spectral radius formula in (3.13), let us consider the case of Arrow–Hurwicz scheme [3], *i.e.* $\theta = 0$. Let $R, J^*$ be locally polyhedral, and $\Sigma_{\overline{L}} = (\sigma_j)_{\{j=1,...,l\}}$ be the singular values of $\overline{L}$, then the eigenvalues of $M_{\text{PD}}$ are

$$\rho_j = \frac{(2 - \gamma_R\gamma_J\sigma_j^2) \pm \sqrt{\gamma_R\gamma_J\sigma_j^2(\gamma_R\gamma_J\sigma_j^2 - 4)}}{2},\ j \in \{1, ..., l\}, \tag{3.14}$$

which apparently are complex ($\gamma_R\gamma_J\sigma_j^2 \leq \gamma_R\gamma_J\|L\|^2 < 1$). Moreover,

$$|\rho_j| = \frac{1}{2}\sqrt{(2 - \gamma_R\gamma_J\sigma_j^2)^2 - \gamma_R\gamma_J\sigma_j^2(\gamma_R\gamma_J\sigma_j^2 - 4)} = 1.$$

This implies that $M_{\text{PD}}$ has multiple eigenvalues with absolute values all equal to 1, then owing to the result of [5], we have $M_{\text{PD}}$ is not convergent.

Furthermore, for $\theta \in [-1, 0[$, we have $1 - \theta\gamma_R\gamma_J\sigma_{\min}^2 > 1$ meaning that $M_{\text{PD}}$ is not convergent, this implies that the correction step proposed in [29] is necessary for $\theta \in [-1, 0]$. Discussion on $\theta > 1$ is left to Section 4.

(iii) If $l \overset{\text{def}}{=} \text{rank}(L) \leq \dim(T_{x^\star}^R)$, then

$$M_{\text{PD}}^\infty = \begin{bmatrix} \mathrm{P}_{S_{x^\star}^R} & \\ & \mathrm{P}_{S_{v^\star}^{J^*}} \end{bmatrix} + \begin{bmatrix} Y & \\ & X \end{bmatrix} \begin{bmatrix} 0_l & & & & \\ & \text{Id}_{p-l} & & & \\ & & 0_{n-p+l} & & \\ & & & \text{Id}_{q-l} & \\ & & & & 0_{m-q} \end{bmatrix} \begin{bmatrix} Y^* & \\ & X^* \end{bmatrix}.$$

**Corollary 3.10.** *Suppose that Algorithm 1 is run under the conditions of Theorem 3.3, then the following claims hold:*

(i) *the linearized iteration (3.10) is equivalent to*

$$(\mathbf{Id} - M_{\mathrm{PD}}^{\infty})(\boldsymbol{z}_{k+1} - \boldsymbol{z}^\star) = M_{\mathrm{PD}}(\mathbf{Id} - M_{\mathrm{PD}}^{\infty})(\boldsymbol{z}_k - \boldsymbol{z}^\star) + o((\mathbf{Id} - M_{\mathrm{PD}}^{\infty})\|\boldsymbol{z}_k - \boldsymbol{z}^\star\|). \tag{3.15}$$

(ii) *If moreover $R, J^*$ are locally polyhedral around the solution pair $(x^\star, v^\star)$, then $M_{\mathrm{PD}}^{\infty}(\boldsymbol{z}_k - \boldsymbol{z}^\star) = 0$ for all $k$ large enough, and (3.15) becomes*

$$\boldsymbol{z}_{k+1} - \boldsymbol{z}^\star = (M_{\mathrm{PD}} - M_{\mathrm{PD}}^{\infty})(\boldsymbol{z}_k - \boldsymbol{z}^\star). \tag{3.16}$$

**Proof.** See [37, Corollary 6.5]. □

### 3.3 Local linear convergence

Finally, we are able to present the local linear convergence of Primal–Dual splitting methods.

**Theorem 3.11 (Local linear convergence).** *Suppose that Algorithm 1 is run under the conditions of Theorem 3.3, then the following holds:*

(i) *given any $\rho \in ]\rho(M_{\mathrm{PD}} - M_{\mathrm{PD}}^{\infty}), 1[$, there exist a $K$ large enough such that for all $k \geq K$,*

$$\|(\mathbf{Id} - M_{\mathrm{PD}}^{\infty})(\boldsymbol{z}_k - \boldsymbol{z}^\star)\| = O(\rho^{k-K}). \tag{3.17}$$

(ii) *If moreover, $R, J^*$ are locally polyhedral around $(x^\star, v^\star)$, then there exists a $K$ large enough such that for all $k \geq K$, we have directly*

$$\|\boldsymbol{z}_k - \boldsymbol{z}^\star\| = O(\rho^{k-K}), \tag{3.18}$$

*for $\rho \in [\rho(M_{\mathrm{PD}} - M_{\mathrm{PD}}^{\infty}), 1[$.*

**Remark 3.12.**

(i) Similar to Proposition 3.6 and Remark 3.7, the above result remains hold if $(\gamma_J, \gamma_R)$ are varying yet convergent. However, the local rate convergence of $\|\boldsymbol{z}_k - \boldsymbol{z}^\star\|$ will depends on how fast $\{(\gamma_{J,k}, \gamma_{R,k})\}_k$ converge, that means, if they converge at a sublinear rate, then the convergence rate of $\|\boldsymbol{z}_k - \boldsymbol{z}^\star\|$ will eventually become sublinear. See [37, Section 8.3] for the case of Douglas–Rachford splitting method.

(ii) When $F = G^* = 0$ and both $R$ and $J^*$ are locally polyhedral around the $(x^\star, v^\star)$, then the convergence rate of the Primal–Dual splitting method is controlled by $\theta$ and $\gamma_J\gamma_R$ as shown in (3.13); see the upcoming section for a detailed discussion.

For general situations (*i.e.* $F, G^*$ are nontrivial and $R, J^*$ are general partly smooth functions), the factors that contribute to the local convergence rate are much more complicated, which involve the Riemannian Hessians of the functions.

# 4 Discussions

In this part, we present several discussions on the obtained local linear convergence result, including acceleration, the effects of $\theta \geq 1$, local oscillation and comparison with FB splitting method and DR splitting method.

To make the result easier to deliver, for the rest of this section we focus on the case where $F = G^* = 0$, *i.e.* the Primal–Dual splitting method of [13], and moreover $R, J^*$ are locally polyhedral around the solution pair $(x^\star, v^\star)$. Under such setting, the matrix defined in (3.9) becomes

$$M_{\mathrm{PD}} \stackrel{\mathrm{def}}{=} \begin{bmatrix} \mathrm{Id}_n & -\gamma_R \bar{L}^* \\ \gamma_J \bar{L} & \mathrm{Id}_m - (1+\theta)\gamma_J\gamma_R \bar{L}\bar{L}^* \end{bmatrix}. \tag{4.1}$$

### 4.1 Choice of $\theta$

Owing to Lemma 3.8, the matrix $M_{\mathrm{PD}}$ in (4.1) is convergent for $\theta \in ]0,1]$, see Eq. (3.12), with the spectral radius

$$\rho(M_{\mathrm{PD}} - M_{\mathrm{PD}}^{\infty}) = \sqrt{1 - \theta\gamma_R\gamma_J\sigma_{\min}^2} < 1, \tag{4.2}$$

with $\sigma_{\min}$ being the smallest non-zero singular value of $\bar{L}$.

In general, given a solution pair $(x^\star, v^\star)$, $\sigma_{\min}$ is fixed, hence the spectral radius $\rho(M_{\mathrm{PD}} - M_{\mathrm{PD}}^{\infty})$ is simply controlled by $\theta$ and the product $\gamma_J\gamma_R$. To make the local convergence rate as faster as possible, it is obvious that we need to make the value of $\theta\gamma_J\gamma_R$ as big as possible. Recall in the global convergence of Primal–Dual splitting method or the result from [13], that

$$\gamma_J\gamma_R\|L\|^2 < 1.$$

Denote $\sigma_{\max}$ the biggest singular value of $\bar{L}$. It is then straightforward that $\gamma_J\gamma_R\sigma_{\max}^2 \leq \gamma_J\gamma_R\|L\|^2 < 1$ and moreover

$$\begin{aligned}\rho(M_{\mathrm{PD}} - M_{\mathrm{PD}}^{\infty}) &= \sqrt{1 - \theta\gamma_R\gamma_J\sigma_{\min}^2} \\ &> \sqrt{1 - \theta(\sigma_{\min}/\|L\|)^2} \geq \sqrt{1 - \theta(\sigma_{\min}/\sigma_{\max})^2}.\end{aligned} \tag{4.3}$$

If we define $\mathrm{cnd} \stackrel{\mathrm{def}}{=} \sigma_{\max}/\sigma_{\min}$ the condition number of $L$, then we have

$$\rho(M_{\mathrm{PD}} - M_{\mathrm{PD}}^{\infty}) > \sqrt{1 - \theta(1/\mathrm{cnd})^2}.$$

To this end, it is clear that $\theta = 1$ gives the best convergence rate for $\theta \in [-1, 1]$. Next let us look at what happens locally if we choose $\theta > 1$.The spectral radius formula (4.2) implies that bigger value of $\theta$ yields smaller spectral radius $\rho(M_{\mathrm{PD}} - M_{\mathrm{PD}}^{\infty})$. Therefore, locally we should choose $\theta$ as big as possible. However, there is an upper bound of $\theta$ which is discussed below.

Following Remark 3.9, let $\Sigma_{\bar{L}} = (\sigma_j)_{\{j=1,\dots,l\}}$ be the singular values of $\bar{L}$, let $\rho_j$ be the eigenvalue of $M_{\mathrm{PD}} - M_{\mathrm{PD}}^{\infty}$, we have known that $\rho_j$ is complex with

$$\rho_j = \frac{\big(2 - (1+\theta)\gamma_R\gamma_J\sigma_j^2\big) \pm \sqrt{(1+\theta)^2\gamma_R^2\gamma_J^2\sigma_j^4 - 4\gamma_R\gamma_J\sigma_j^2}}{2}, \quad |\rho_j| = \sqrt{1 - \theta\gamma_R\gamma_J\sigma_j^2}.$$

Now let $\theta > 1$, then in order to let $|\rho_j|$ make sense for all $j \in \{1, \dots, l\}$, there must holds

$$1 - \theta\gamma_R\gamma_J\sigma_{\max}^2 \geq 0 \Longleftrightarrow \theta \leq \frac{1}{\gamma_R\gamma_J\sigma_{\max}^2},$$

which means that $\theta$ indeed is bounded from above.

Unfortunately, since $\bar{L} = \mathrm{P}_{T_{x^\star}^R} L \mathrm{P}_{T_{v^\star}^{J^*}}$, the upper bound can be only obtained if we had the solution pair $(x^\star, v^\star)$. However, in practice we can use back-tracking or the Armijo-Goldstein-rule to find the proper $\theta$. See Section 6.4 for an illustration of online searching of $\theta$. It should be noted that such updating rule can also be applied to $\gamma_J, \gamma_R$ since we have $\|\bar{L}\| \leq \|L\|$, hence locally bigger values of them should be applicable. Moreover, it should be noted that in practice one can choose to enlarge either $\theta$ or $\gamma_J\gamma_R$ as they will have very similar acceleration outcome.

**Remark 4.1.** It should be noted that the above discussion on the effect of $\theta > 1$ may only valid for the case $F = 0, G^* = 0$, *i.e.* the Primal–Dual splitting method of [13]. If $F$ and/or $G^*$ are not vanished, then locally, $\theta < 1$ may give faster convergence rate.

## 4.2 Oscillations

For the inertial FB splitting and FISTA [7] methods, it is shown in [36] that they locally oscillate when the inertia momentum are too high (see [36, Section 4.4] for more details). When solving certain type of problems (*i.e.* $F = G^* = 0$ and $R, J^*$ are locally polyhedral around the solution pair $(x^\star, v^\star)$), the Primal–Dual splitting method also locally oscillates (see Figure 6 for an illustration). As revealed in the proof of Lemma 3.8, all the eigenvalues of $M_{\mathrm{PD}} - M_{\mathrm{PD}}^\infty$ in (4.1) are complex. This means that locally the sequences generated by the Primal–Dual splitting iteration may oscillate.

For $\sigma_{\min}$, the smallest non-zero singular of $\bar{L}$, one of its corresponding eigenvalues of $M_{\mathrm{PD}}$ reads

$$\rho_{\sigma_{\min}} = \frac{(2-(1+\theta)\gamma_J\gamma_R\sigma_{\min}^2) + \sqrt{(1+\theta)^2\gamma_R^2\gamma_J^2\sigma_{\min}^4 - 4\gamma_J\gamma_R\sigma_{\min}^2}}{2},$$

and $(1+\theta)^2\gamma_R^2\gamma_J^2\sigma_{\min}^4 - 4\gamma_J\gamma_R\sigma_{\min}^2 < 0$. Denote $\omega$ the argument of $\rho_{\sigma_{\min}}$, then

$$\cos(\omega) = \frac{2-(1+\theta)\gamma_J\gamma_R\sigma_{\min}^2}{\sqrt{1-\theta\gamma_J\gamma_R\sigma_{\min}^2}}. \tag{4.4}$$

The oscillation period of the sequence $\|\boldsymbol{z}_k - \boldsymbol{z}^\star\|$ is then exactly $\frac{\pi}{\omega}$. See Figure 6 for an illustration.

**Remark 4.2.**

(i) Complex eigenvalues is only a necessary condition for local oscillation behaviour. When the involved functions are all polyhedral, the eigenvalues of the local linearized operator of DR splitting method are also complex (see [24, 4]), however, the iterates of DR splitting method do not oscillate.

(ii) The mechanisms of oscillation between Primal–Dual splitting and the inertial FB splitting (including FISTA) methods are different. The oscillation of FISTA is caused by the inertial momentum being too large, while the oscillation of Primal–Dual splitting is due to the polyhedrality of the functions, which makes the eigenvalues all complex (see Lemma 3.8). Furthermore, if $F \neq 0$ and/or $G^* \neq 0$, then the Primal–Dual splitting method may not oscillate, see Figure 2 (c) numerical evidence.

## 4.3 Relations with FB and DR/ADMM

In this part, we discuss the relation between the obtained result and our previous work on local linear convergence of FB splitting [34, 36] and DR splitting [38, 37].

### 4.3.1 FB splitting

For problem ($\mathcal{P}_{\mathrm{P}}$), when $J = G^* = 0$, Algorithm 1 reduces to, denoting $\gamma = \gamma_R$ and $\beta = \beta_F$,

$$x_{k+1} = \mathrm{prox}_{\gamma R}\big(x_k - \gamma\nabla F(x_k)\big),\ \gamma \in ]0, 2\beta[, \tag{4.5}$$

which is the non-relaxed FB splitting method [39] with constant step-size.

Let $x^\star \in \mathrm{Argmin}(R+F)$ be a global minimizer to which $\{x_k\}_{k\in\mathbb{N}}$ of (4.5) converges, the non-degeneracy condition (ND) for identification then becomes

$$-\nabla F(x^\star) \in \mathrm{ri}\big(\partial R(x^\star)\big),$$

which recovers the conditions of [36, Theorem 4.11]. Following the notations of Section 3, define $M_{\mathrm{FB}} \stackrel{\mathrm{def}}{=} W_{\bar{R}}(\mathrm{Id}_n - \gamma H_F)$, we have for all $k$ large enough

$$x_{k+1} - x^\star = M_{\mathrm{FB}}(x_k - x^\star) + o(\|x_k - x^\star\|).$$

From Theorem 3.11, we obtain the following result for the FB splitting method, for the case $\gamma$ being fixed.

**Corollary 4.3.** *For problem ($\mathcal{P}_{\mathrm{P}}$), let $J = G^* = 0$ and suppose that (**A.1**) holds and* $\mathrm{Argmin}(R + F) \neq \emptyset$*, and the FB iteration (4.5) creates a sequence $x_k \to x^\star \in \mathrm{Argmin}(R + F)$ such that $R \in \mathrm{PSF}_{x^\star}(\mathcal{M}_{x^\star})$, $F$ is $C^2$ near $x^\star$, and condition $-\nabla F(x^\star) \in \mathrm{ri}(\partial R(x^\star))$ holds. Then*

(i) *given any $\rho \in ]\rho(M_{\mathrm{FB}} - M_{\mathrm{FB}}^\infty), 1[$, there exist a $K$ large enough such that for all $k \geq K$,*

$$\|(\mathrm{Id} - M_{\mathrm{FB}}^\infty)(x_k - x^\star)\| = O(\rho^{k-K}). \tag{4.6}$$

(ii) *If moreover, $R$ are locally polyhedral around $x^\star$, there exist a $K$ large enough such that for all $k \geq K$, we have directly*

$$\|x_k - x^\star\| = O(\rho^{k-K}), \tag{4.7}$$

*for $\rho \in [\rho(M_{\mathrm{FB}} - M_{\mathrm{FB}}^\infty), 1[$.*

**Proof.** Owing to [6], $M_{\mathrm{FB}}$ is $\frac{2\beta}{4\beta-\gamma}$-averaged non-expansive, hence convergent. The convergence rates in (4.6) and (4.7) are straightforward from Theorem 3.11. □

Unlike the result in [34, 36], for the local linear convergence of FB splitting method, a so-called restricted injectivity condition (RI) is required, which means that $H_F$ should be positive definite. Moreover, RI condition can be only removed when $J$ is locally polyhedral around $x^\star$ (*e.g.* see [36, Theorem 4.9]), while in this paper, we show that no polyhedrality is needed. However, the price of removing such a condition is that the obtained convergence rate is on a different criterion (*i.e.* $\|(\mathrm{Id} - M_{\mathrm{FB}}^\infty)(x_k - x^\star)\|$) other than the sequence itself (*i.e.* $\|x_k - x^\star\|$).

### 4.3.2 DR splitting and ADMM

Let $F = G^* = 0$ and $L = \mathrm{Id}$, then problem ($\mathcal{P}_{\mathrm{P}}$) becomes

$$\min_{x \in \mathbb{R}^n} R(x) + J(x).$$

For the above problem, below we briefly show that DR splitting is the limiting case of Primal–Dual splitting by letting $\gamma_R \gamma_J = 1$. First, for the Primal–Dual splitting scheme of Algorithm 1, let $\theta = 1$ and change the order of updating the variables, we obtain the following iteration

$$\left\lfloor \begin{aligned} v_{k+1} &= \mathrm{prox}_{\gamma_J J^*}(v_k + \gamma_J \bar{x}_k) \\ x_{k+1} &= \mathrm{prox}_{\gamma_R R}(x_k - \gamma_R v_{k+1}) \\ \bar{x}_{k+1} &= 2x_{k+1} - x_k. \end{aligned} \right. \tag{4.8}$$

Then apply the Moreau's identity (2.4) to $\mathrm{prox}_{\gamma_J J^*}$, let $\gamma_J = 1/\gamma_R$ and define $z_{k+1} = x_k - \gamma_R v_{k+1}$, iteration (4.8) becomes

$$\left\lfloor \begin{aligned} u_{k+1} &= \mathrm{prox}_{\gamma_R J}(2x_k - z_k) \\ z_{k+1} &= z_k + u_{k+1} - x_k \\ x_{k+1} &= \mathrm{prox}_{\gamma_R R}(z_{k+1}), \end{aligned} \right. \tag{4.9}$$

which is the non-relaxed DR splitting method [25]. At convergence, we have $u_k, x_k \to x^\star = \mathrm{prox}_{\gamma_R R}(z^\star)$ where $z^\star$ is a fixed point of the iteration. See also the discussions in [13, Section 4.2].

Specializing the derivation of (3.9) to (4.8) and (4.9), we obtain the following two linearized fixed-point operator for (4.8) and (4.9) respectively

$$M_1 = \begin{bmatrix} \mathrm{Id}_n & -\gamma_R \mathrm{P}_{T_{x^\star}^R} \mathrm{P}_{T_{v^\star}^{J^*}} \\ \gamma_J \mathrm{P}_{T_{v^\star}^{J^*}} \mathrm{P}_{T_{x^\star}^R} & \mathrm{Id}_n - 2\gamma_J \gamma_R \mathrm{P}_{T_{v^\star}^{J^*}} \mathrm{P}_{T_{x^\star}^R} \mathrm{P}_{T_{v^\star}^{J^*}} \end{bmatrix},$$

$$M_2 = \begin{bmatrix} \mathrm{Id}_n & -\gamma_R \mathrm{P}_{T_{x^\star}^R} \mathrm{P}_{T_{v^\star}^{J^*}} \\ \frac{1}{\gamma_R} \mathrm{P}_{T_{v^\star}^{J^*}} \mathrm{P}_{T_{x^\star}^R} & \mathrm{Id}_n - 2\mathrm{P}_{T_{v^\star}^{J^*}} \mathrm{P}_{T_{x^\star}^R} \mathrm{P}_{T_{v^\star}^{J^*}} \end{bmatrix}.$$

Owing to (ii) of Lemma 3.8, $M_1, M_2$ are convergent. Let $\omega$ be the largest principal angle (yet smaller than $\pi/2$) between tangent spaces $T^R_{x^\star}$ and $T^{J^*}_{v^\star}$, then we have the spectral radius of $M_1 - M_1^\infty$ reads ((i) of Remark 3.9),

$$\begin{aligned}\rho(M_1 - M_1^\infty) &= \sqrt{1 - \gamma_J \gamma_R \cos^2(\omega)} \\ &\geq \sqrt{1 - \cos^2(\omega)} = \sin(\omega) = \cos(\pi/2 - \omega).\end{aligned} \tag{4.10}$$

Suppose that the Kuhn-Tucker pair $(x^\star, v^\star)$ is unique, and moreover that $R$ and $J$ are polyhedral. Therefore, we have that if $J^*$ is locally polyhedral near $v^\star$ along $v^\star + T^{J^*}_{v^\star}$, then $J$ is locally polyhedral near $x^\star$ around $x^\star + T^J_{x^\star}$, and moreover there holds $T^J_{x^\star} = (T^{J^*}_{v^\star})^\perp$. As a result, the principal angles between $T^R_{x^\star}, T^{J^*}_{v^\star}$ and the ones between $T^R_{x^\star}, T^J_{x^\star}$ are complementary, which means that $\pi/2 - \omega$ is the Friedrichs angle between tangent spaces $T^R_{x^\star}, T^J_{x^\star}$. Thus, following (4.10), we have

$$\rho(M_1 - M_1^\infty) = \sqrt{1 - \gamma_J \gamma_R \cos^2(\omega)} \geq \cos(\pi/2 - \omega) = \rho(M_2 - M_2^\infty).$$

# 5 Multiple infimal convolutions

In this section, we consider problem $(\mathcal{P}_\mathrm{P})$ with more than one infimal convolution. Le $m \geq 1$ be a positive integer. Consider the problem of minimizing

$$\min_{x\in\mathbb{R}^n} R(x) + F(x) + \textstyle\sum_{i=1}^m (J_i \mathbin{\triangledown} G_i)(L_i x), \tag{$\mathcal{P}^m_\mathrm{P}$}$$

where (A.1) holds for $R$ and $F$, and for every $i = 1, ..., m$ the followings are hold:

**(A'.2)** $J_i, G_i \in \Gamma_0(\mathbb{R}^{m_i})$, with $G_i$ being differentiable and $\beta_{G,i}$-strongly convex for $\beta_{G,i} > 0$.

**(A'.3)** $L_i : \mathbb{R}^n \to \mathbb{R}^{m_i}$ is a linear operator.

**(A'.4)** The condition

$$0 \in \mathrm{ran}\left(\partial R + \nabla F + \textstyle\sum_{i=1}^m L_i^*(\partial J_i \Box \partial G_i) L_i\right).$$

The dual problem of $(\mathcal{P}^m_\mathrm{P})$ reads,

$$\min_{v_1\in\mathbb{R}^{m_1},...,v_m\in\mathbb{R}^{m_m}} (R^* \mathbin{\triangledown} F^*)\big(-\textstyle\sum_{i=1}^m L_i^* v_i\big) + \textstyle\sum_{i=1}^m \big(J_i^*(v_i) + G_i^*(v_i)\big). \tag{$\mathcal{P}^m_\mathrm{D}$}$$

Problem $(\mathcal{P}^m_\mathrm{P})$ is considered in [50, 19], and a Primal–Dual splitting algorithm is proposed there which is an extension of Algorithm 1 using a product space trick, see Algorithm 2 hereafter for details. In both schemes proposed by [50, 19], the choice of $\theta$ is set as 1.

**Algorithm 2:** A Primal–Dual splitting method

**Initial**: Choose $\gamma_R, (\gamma_{J_i})_i > 0$. For $k = 0$, $x_0 \in \mathbb{R}^n$, $v_{i,0} \in \mathbb{R}^{m_i}$, $i \in \{1, ..., m\}$;

**repeat**

$$\left\lfloor \begin{aligned} & x_{k+1} = \mathrm{prox}_{\gamma_R R}\big(x_k - \gamma_R \nabla F(x_k) - \gamma_R \textstyle\sum_i L_i^* v_{i,k}\big) \\ & \bar{x}_{k+1} = 2x_{k+1} - x_k \\ & \text{For } i = 1, ..., m \\ & \left\lfloor v_{i,k+1} = \mathrm{prox}_{\gamma_{J_i} J_i^*}\big(v_{i,k} - \gamma_{J_i} \nabla G_i^*(v_{i,k}) + \gamma_{J_i} L_i \bar{x}_{k+1}\big), \right. \end{aligned} \right. \tag{5.1}$$

$k = k + 1$;

**until** *convergence*;

## 5.1 Product space

The following result is taken from [19]. Define the product space $\boldsymbol{\mathcal{K}} = \mathbb{R}^n \times \mathbb{R}^{m_1} \times \cdots \times \mathbb{R}^{m_m}$, and let $\mathbf{Id}$ be the identity operator on $\boldsymbol{\mathcal{K}}$. Define the following operators

$$\boldsymbol{A} \stackrel{\text{def}}{=} \begin{bmatrix} \partial R & L_1^* & \cdots & L_m^* \\ -L_1 & \partial J_1 & & \\ \vdots & & \ddots & \\ -L_m & & & \partial J_m \end{bmatrix}, \boldsymbol{B} \stackrel{\text{def}}{=} \begin{bmatrix} \nabla F & & & \\ & \nabla G_1^* & & \\ & & \ddots & \\ & & & \nabla G_m^* \end{bmatrix}, \boldsymbol{\mathcal{V}} \stackrel{\text{def}}{=} \begin{bmatrix} \frac{\mathrm{Id}_n}{\gamma_R} & -L_1^* & \cdots & -L_m^* \\ -L_1 & \frac{\mathrm{Id}_{m_1}}{\gamma_{J_1}} & & \\ \vdots & & \ddots & \\ -L_m & & & \frac{\mathrm{Id}_{m_m}}{\gamma_{J_m}} \end{bmatrix}. \tag{5.2}$$

Then $\boldsymbol{A}$ is maximal monotone, $\boldsymbol{B}$ is $\min\{\beta_F, \beta_{G_1}, ..., \beta_{G_m}\}$-cocoercive, and $\boldsymbol{\mathcal{V}}$ is symmetric and $\nu$-positive definite with $\nu = (1-\sqrt{\gamma_R \sum_i \gamma_{J_i} \|L_i\|^2}) \min\{\frac{1}{\gamma_R}, \frac{1}{\gamma_{J_1}}, ..., \frac{1}{\gamma_{J_m}}\}$. Define $\boldsymbol{z}_k = (x_k, v_{1,k}, \cdots, v_{m,k})^T$, then it can be shown that (5.1) is equivalent to

$$\boldsymbol{z}_{k+1} = (\boldsymbol{\mathcal{V}} + \boldsymbol{A})^{-1}(\boldsymbol{\mathcal{V}} - \boldsymbol{B})\boldsymbol{z}_k = (\mathbf{Id} + \boldsymbol{\mathcal{V}}^{-1}\boldsymbol{A})^{-1}(\mathbf{Id} - \boldsymbol{\mathcal{V}}^{-1}\boldsymbol{B})\boldsymbol{z}_k. \tag{5.3}$$

## 5.2 Local convergence analysis

Let $(x^\star, v_1^\star, ..., v_m^\star)$ be a Kuhn-Tucker pair. Define the following functions

$$\overline{J_i^*}(v_i) \stackrel{\text{def}}{=} J_i^*(v_i) - \langle v_i,\, L_i x^\star - \nabla G_i^*(v_i^\star)\rangle,\ \ v_i \in \mathbb{R}^{m_i},\ \ i \in \{1, ..., m\}, \tag{5.4}$$

and the Riemannian Hessian of each $\overline{J_i^*}$,

$$H_{\overline{J_i^*}} \stackrel{\text{def}}{=} \gamma_{J_i} \mathrm{P}_{T_{v_i^\star}^{J_i^*}} \nabla^2_{\mathcal{M}_{v_i^\star}^{J_i^*}} \overline{J_i^*}(v_i^\star) \mathrm{P}_{T_{v_i^\star}^{J_i^*}} \ \text{ and } \ W_{\overline{J_i^*}} \stackrel{\text{def}}{=} (\mathrm{Id}_{m_i} + H_{\overline{J_i^*}})^{-1},\ \ i \in \{1, ..., m\}. \tag{5.5}$$

For each $i \in \{1, ..., m\}$, owing to Lemma 3.5, we have that $W_{\overline{J_i^*}}$ is firmly non-expansive if the non-degeneracy condition ($\mathrm{ND}_m$) below holds. Now suppose that $F$ locally is $C^2$ around $x^\star$ and $G_i^*$ locally is $C^2$ around $v_i^\star$, define the restricted Hessian $H_{G_i^*} \stackrel{\text{def}}{=} \mathrm{P}_{T_{v_i^\star}^{J_i^*}} \nabla^2 G_i^*(v_i^\star) \mathrm{P}_{T_{v_i^\star}^{J_i^*}}$. Define $\overline{H}_{G_i^*} \stackrel{\text{def}}{=} \mathrm{Id}_{m_i} - \gamma_{J_i^*} H_{G_i^*}$, $\overline{L}_i \stackrel{\text{def}}{=} \mathrm{P}_{T_{v_i^\star}^{J_i^*}} L_i \mathrm{P}_{T_{x^\star}^R}$, and the matrix

$$M_{\mathrm{PD}} \stackrel{\text{def}}{=} \begin{bmatrix} W_{\overline{R}}\overline{H}_F & -\gamma_R W_{\overline{R}}\overline{L}_1^* & \cdots & -\gamma_R W_{\overline{R}}\overline{L}_m^* \\ \gamma_{J_1^*} W_{J_1^*}\overline{L}_1(2W_{\overline{R}}\overline{H}_F - \mathrm{Id}_n) & W_{J_1^*}(\overline{H}_{G_1^*} - 2\gamma_{J_1^*}\gamma_R \overline{L}_1 W_{\overline{R}}\overline{L}_1^*) & & \\ \vdots & & \ddots & \\ \gamma_{J_m^*} W_{J_m^*}\overline{L}_m(2W_{\overline{R}}\overline{H}_F - \mathrm{Id}_n) & \cdots & & W_{J_m^*}(\overline{H}_{G_m^*} - 2\gamma_{J_m^*}\gamma_R \overline{L}_m W_{\overline{R}}\overline{L}_m^*) \end{bmatrix}. \tag{5.6}$$

Using the same strategy of the proof of Lemma 3.8, one can show that $M_{\mathrm{PD}}$ is convergent, which again is denoted as $M_{\mathrm{PD}}^\infty$, and $\rho(M_{\mathrm{PD}} - M_{\mathrm{PD}}^\infty) < 1$.

**Corollary 5.1.** *Consider Algorithm 2 under assumptions (A.1) and (A'.2)-(A'.4). Choose $\gamma_R, (\gamma_{J_i})_i > 0$ such that*

$$2\min\{\beta_F, \beta_{G_1}, ..., \beta_{G_m}\} \min\Big\{\frac{1}{\gamma_R}, \frac{1}{\gamma_{J_1}}, ..., \frac{1}{\gamma_{J_m}}\Big\}\Big(1 - \sqrt{\gamma_R \textstyle\sum_i \gamma_{J_i}\|L_i\|^2}\Big) > 1. \tag{5.7}$$

*$(x_k, v_{1,k}, ..., v_{m,k}) \to (x^\star, v_1^\star, ..., v_m^\star)$, where $x^\star$ solves ($\mathcal{P}_{\mathrm{P}}^m$) and $(v_1^\star, ..., v_m^\star)$ solve ($\mathcal{P}_{\mathrm{D}}^m$). If moreover $R \in \mathrm{PSF}_{x^\star}(\mathcal{M}_{x^\star}^R)$ and $J_i^* \in \mathrm{PSF}_{v_i^\star}(\mathcal{M}_{v_i^\star}^{J_i^*})$, $i \in \{1, ..., m\}$, and the non-degeneracy condition holds*

$$\begin{aligned} -\textstyle\sum_i L_i^* v_i^\star - \nabla F(x^\star) &\in \mathrm{ri}\big(\partial R(x^\star)\big) \\ L_i x^\star - \nabla G_i^*(v^\star) &\in \mathrm{ri}\big(\partial J_i^*(v^\star)\big),\ \ \forall i \in \{1, ..., m\}. \end{aligned} \tag{$\mathrm{ND}_m$}$$

*Then,*

(i) $\exists K > 0$ *such that for all* $k \geq K$,

$$(x_k, v_{1,k}, ..., v_{m,k}) \in \mathcal{M}^R_{x^\star} \times \mathcal{M}^{J_1^*}_{v_1^\star} \times \cdots \times \mathcal{M}^{J_m^*}_{v_m^\star}.$$

(ii) *Given any* $\rho \in ]\rho(M_{\mathrm{PD}} - M^\infty_{\mathrm{PD}}), 1[$, *there exist a* $K$ *large enough such that for all* $k \geq K$,

$$\|(\mathbf{Id} - M^\infty_{\mathrm{PD}})(\boldsymbol{z}_k - \boldsymbol{z}^\star)\| = O(\rho^{k-K}). \tag{5.8}$$

*If moreover,* $R, J_1^*, ..., J_m^*$ *are locally polyhedral around* $(x^\star, v_1^\star, ..., v_m^\star)$, *then we have directly have* $\|\boldsymbol{z}_k - \boldsymbol{z}^\star\| = O(\rho^{k-K})$.

# 6 Numerical experiments

In this section, we illustrate our theoretical results on several concrete examples arising from fields including inverse problem, signal/image processing and machine learning.

## 6.1 Examples of partly smooth function

Table 1 provides some examples of partly smooth functions that we will use throughout this section, more details about them can be found in [36, Section 5] and the references therein.

Table 1: Examples of partly smooth functions. For $x \in \mathbb{R}^n$ and some subset of indices $\mathfrak{b} \subset \{1, \ldots, n\}$, $x_\mathfrak{b}$ is the restriction of $x$ to the entries indexed in $\mathfrak{b}$. $D_{\mathrm{DIF}}$ stands for the finite differences operator.

| Function | Expression | Partial smooth manifold |
|---|---|---|
| $\ell_1$-norm | $\|x\|_1 = \sum_{i=1}^n \|x_i\|$ | $\mathcal{M} = T_x = \{z \in \mathbb{R}^n : I_z \subseteq I_x\}$, $I_x = \{i : x_i \neq 0\}$ |
| $\ell_{1,2}$-norm | $\sum_{i=1}^m \|x_{\mathfrak{b}_i}\|$ | $\mathcal{M} = T_x = \{z \in \mathbb{R}^n : I_z \subseteq I_x\}$, $I_x = \{i : x_{\mathfrak{b}_i} \neq 0\}$ |
| $\ell_\infty$-norm | $\max_{i=\{1,\ldots,n\}} \|x_i\|$ | $\mathcal{M} = T_x = \{z \in \mathbb{R}^n : z_{I_x} \in \mathbb{R}\mathrm{sign}(x_{I_x})\}$, $I_x = \{i : \|x_i\| = \|x\|_\infty\}$ |
| TV semi-norm | $\|x\|_{\mathrm{TV}} = \|D_{\mathrm{DIF}} x\|_1$ | $\mathcal{M} = T_x = \{z \in \mathbb{R}^n : I_{D_{\mathrm{DIF}} z} \subseteq I_{D_{\mathrm{DIF}} x}\}$, $I_{D_{\mathrm{DIF}} x} = \{i : (D_{\mathrm{DIF}} x)_i \neq 0\}$ |
| Nuclear norm | $\|x\|_* = \sum_{i=1}^r \sigma(x)$ | $\mathcal{M} = \{z \in \mathbb{R}^{n_1 \times n_2} : \mathrm{rank}(z) = \mathrm{rank}(x) = r\}$, $\sigma(x)$ singular values of $x$ |

The $\ell_1$, $\ell_\infty$-norms and the anisotropic TV semi-norm are all polyhedral functions, hence their Riemannian Hessian are simply 0. The $\ell_{1,2}$-norm is not polyhedral yet partly smooth relative to a subspace, the nuclear norm is partly smooth relative to the set of fixed-rank matrices, which on the other hand is curved, the Riemannian of these two functions are non-trivial and can be found in [49] and references therein.

## 6.2 Linear inverse problems

Given an object $x_{\mathrm{ob}} \in \mathbb{R}^n$, often times we can not access it directly, but through the following observation model,

$$b = \mathcal{K} x_{\mathrm{ob}}, \tag{6.1}$$

where $b \in \mathbb{R}^m$ is the observation, $\mathcal{K} : \mathbb{R}^n \to \mathbb{R}^m$ is some linear operator. A more complicated situation is when the observation is contaminated by noise, namely, $b = \mathcal{K} x_{\mathrm{ob}} + w$, where $w \in \mathbb{R}^m$ is the noise.

The operator $\mathcal{K}$ usually is ill-conditioned or even singular, hence recovering or approximating $x_{\mathrm{ob}}$ from (6.1) in general is ill-posed. However, usually some prior knowledge of $x_{\mathrm{ob}}$ can be available. Thus, a popular approach to recover $x_{\mathrm{ob}}$ from $b$ is via regularization, by solving

$$\min_{x \in \mathbb{R}^n} R(x) + J(\mathcal{K}x - b), \tag{6.2}$$

where

- $R \in \Gamma_0(\mathbb{R}^n)$ is the regularizer based on the prior information, *e.g.* $\ell_1, \ell_{1,2}, \ell_\infty$-norms, nuclear norm;
- $J \in \Gamma_0(\mathbb{R}^m)$ enforces fidelity to the observations. Typically $J = \iota_0$ when there is no noise, *i.e.* $w = 0$.

Clearly, (6.2) is a special instance of ($\mathcal{P}_{\mathrm{P}}$) with $F = G^* = 0$, hence Algorithm 1 can be applied.

We consider problem (6.2) with $R$ being $\ell_1, \ell_{1,2}, \ell_\infty$-norms, and nuclear norm. $\mathcal{K} \in \mathbb{R}^{m\times n}$ is generated uniformly at random with independent zero-mean standard Gaussian entries. The settings of the experiments are:

**$\ell_1$-norm** $(m, n) = (48, 128)$, $\|x_{\mathrm{ob}}\|_0 = 8$;
**$\ell_{1,2}$-norm** $(m, n) = (48, 128)$, $x_{\mathrm{ob}}$ has 3 non-zero blocks of size 4;
**$\ell_\infty$-norm** $(m, n) = (63, 64)$, $|I(x_{\mathrm{ob}})| = 8$;
**Nuclear norm** $(m, n) = (500, 1024)$, $x_{\mathrm{ob}} \in \mathbb{R}^{32\times 32}$ and $\mathrm{rank}(x_{\mathrm{ob}}) = 4$.

Figure 1 displays the profile of $\|\boldsymbol{z}_k - \boldsymbol{z}^\star\|$ as a function of $k$, and the starting point of the dashed line is the iteration number at which the active partial smoothness manifold of $\mathcal{M}^R_{x^\star}$ is identified (recall that $\mathcal{M}^{J^*}_{v^\star} = \{0\}$ which is trivially identified from the first iteration). One can easily see that for the $\ell_1$ and $\ell_\infty$ norms, Theorem 3.11 applies and our estimates are very tight, meaning that the dashed and solid lines has the same slope. For the case of $\ell_{1,2}$-norm and nuclear norm, though not optimal, our estimates are very tight.

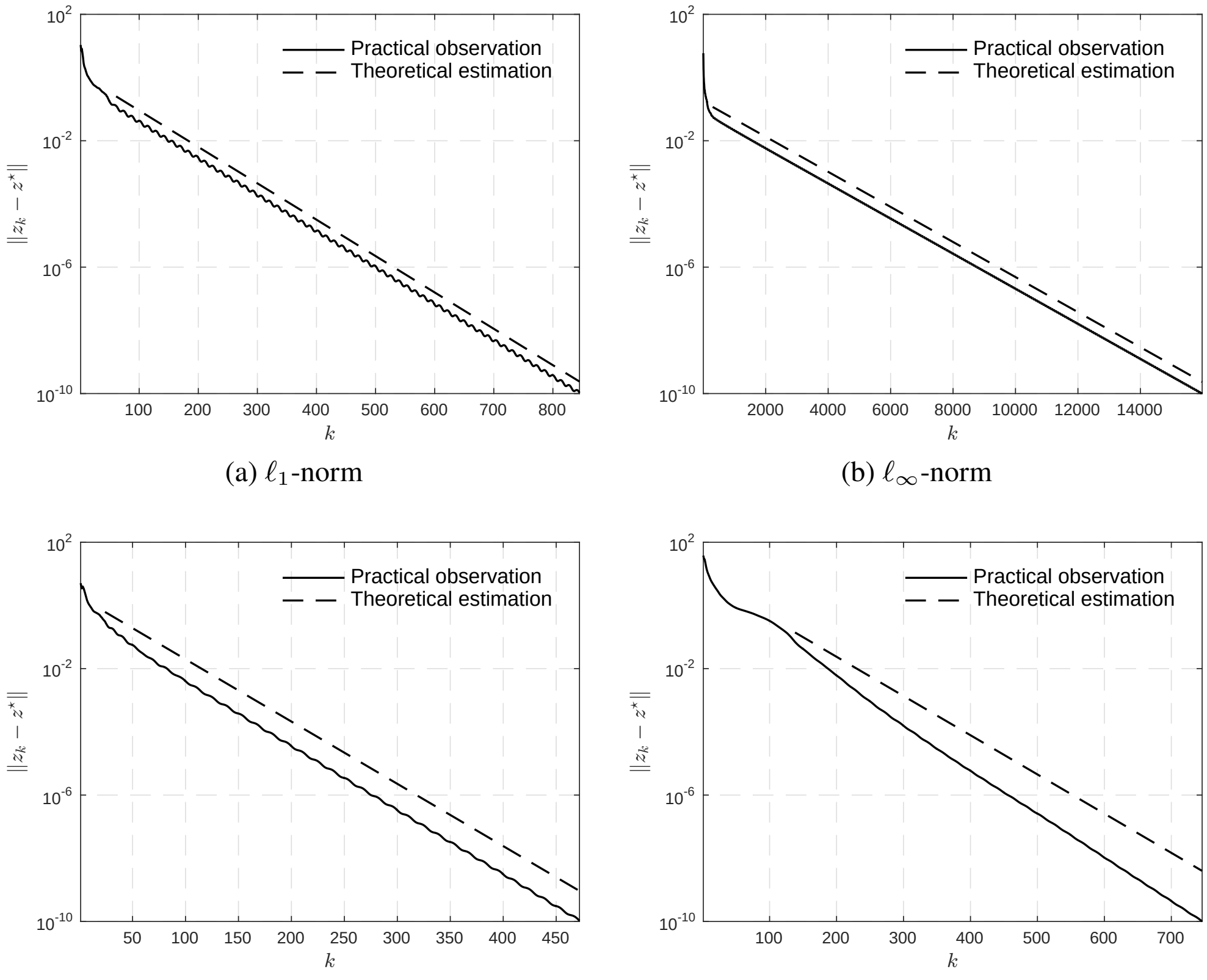

(a) $\ell_1$-norm (b) $\ell_\infty$-norm (c) $\ell_{1,2}$-norm (d) Nuclear norm

Figure 1: Observed (solid) and predicted (dashed) convergence profiles of Algorithm 1 in terms of $\|\boldsymbol{z}_k - \boldsymbol{z}^\star\|$. (a) $\ell_1$-norm. (b) $\ell_\infty$-norm. (c) $\ell_{1,2}$-norm. (d) Nuclear norm. The starting point of the dashed line is the iteration at which the active manifold of $J$ is identified.

### 6.3 Total variation based denoising

In this part, we consider several examples of total variation based denoising, for the first two examples, we suppose that we observe $b = x_{\text{ob}} + w$, where $x_{\text{ob}}$ is a piecewise-constant vector, and $w$ is an unknown noise supposed to be either uniform or sparse. The goal is to recover $x_{\text{ob}}$ from $y$ using the prior information on $x_{\text{ob}}$ (*i.e.* piecewise-smooth) and $w$ (uniform or sparse). To achieve this goal, a popular and natural approach in the signal processing literature is to solve

$$\min_{x \in \mathbb{R}^n} \|D_{\text{DIF}} x\|_1 \quad \text{subject to} \quad \|b - x\|_p \leq \tau, \tag{6.3}$$

where $p = +\infty$ for uniform noise, and $p = 1$ for sparse noise, and $\tau > 0$ is a parameter depending on the noise level.

Problem (6.3) can also formulated into the form of ($\mathcal{P}_{\text{P}}$). Indeed, one can take $R = \iota_{\|b - \cdot\|_p \leq \tau}$, $J = \|\cdot\|_1$, $F = G^* = 0$, and $L = D_{\text{DIF}}$ is the finite difference operator (with appropriate boundary conditions). The proximity operators of $R$ and $J$ can be computed easily. Clearly, both two indicator functions are polyhedral, and their proximal operator are simple to compute.

For both examples, we set $n = 128$ and $x_{\text{ob}}$ is such that $D_{\text{DIF}} x_{\text{ob}}$ has 8 nonzero entries. For $p = +\infty$, $w$ is generated uniformly in $[-1, 1]$, and for $p = 1$, $w$ is sparse with 16 nonzero entries. The corresponding local convergence profiles are depicted in Figure 2(a)-(b). Owing to polyhedrality, our rate predictions are again optimal.

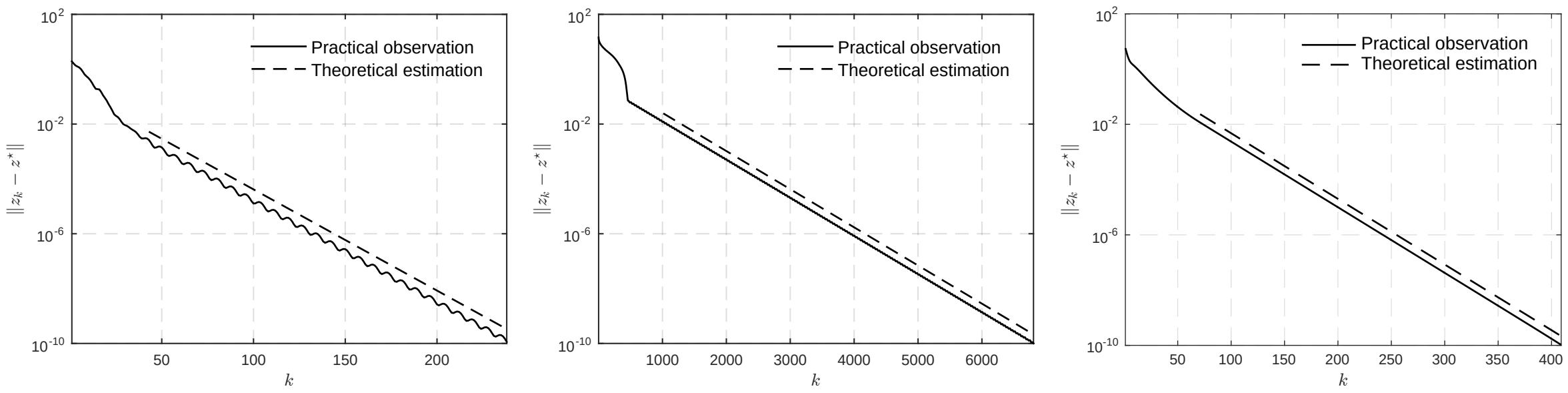

(a) Sparse noise removal $p = 1$ (b) Uniform noise removal $p = \infty$ (c) Gaussian noise removal

Figure 2: Observed (solid) and predicted (dashed) convergence profiles of Primal–Dual (1.2) in terms of $\|\boldsymbol{z}_k - \boldsymbol{z}^\star\|$. (a) Sparse noise removal. (b) Uniform noise removal. (c) Gaussian noise removal. The starting point of the dashed line is the iteration at which the active manifold of $J$ is identified.

We also consider an underdetermined linear regression problem

$$b = \mathcal{K} x_{\text{ob}} + w.$$

We assume that the vectors to promote are group sparse and each non-zero group is constant. This regression problem can then be approached by solving

$$\min_{x \in \mathbb{R}^n} \mu_1 \|x\|_{1,2} + \frac{1}{2}\|\mathcal{K}x - b\|^2 + \mu_2 \|D_{\text{DIF}} x\|_1,$$

where $\mu_i > 0$, $\|\cdot\|_{1,2}$ is a first regularizer that favours group sparsity, and $\|D_{\text{DIF}} \cdot\|_1$ a second regularizer designed to promote piece-wise constancy. This is again in the form of ($\mathcal{P}_{\text{P}}$), where $R = \mu_1\|\cdot\|_{1,2}$, $F = \frac{1}{2}\|\mathcal{K} \cdot - b\|^2$, $J = \mu_2\|\cdot\|_1$, $G^* = 0$, and $L = D_{\text{DIF}}$. For this example, we set $x_{\text{ob}} \in \mathbb{R}^{128}$ with 2 piecewise constant non-zeros blocks of size 8. The result is shown in Figure 2(c), the estimate is not as tight as the other 2 examples, but still sharp enough.

### 6.4 Choices of $\theta$ and $\gamma_J, \gamma_R$

In this part, we present a comparison on different choices of $\theta$ and $\gamma_J, \gamma_R$ to see their influences on the finite identification and local linear convergence rate. Two examples are consider for these comparisons, problem (6.2) with $R$ being $\ell_1$-norm and $\ell_{1,2}$-norm.

**Fixed $\theta$** We consider first the case of fixing $\theta$, and changing the value of $\gamma_J \gamma_R \|L\|^2$. 4 different cases are considered, which are

$$\gamma_J \gamma_R \|L\|^2 \in \{0.3, 0.6, 0.8, 0.99\},$$

and we fix $\theta = 1$, moreover we set $\gamma_J = \gamma_R$. The comparison result is shown in Figure 3, and we have the following observations

(i) The smaller the value of $\gamma_J \gamma_R \|L\|^2$, the slower the iteration converges;

(ii) Bigger value of $\gamma_R$ leads to faster identification (since $J^*$ is globally $C^2$-smooth, so only the identification of $R$ for this case).

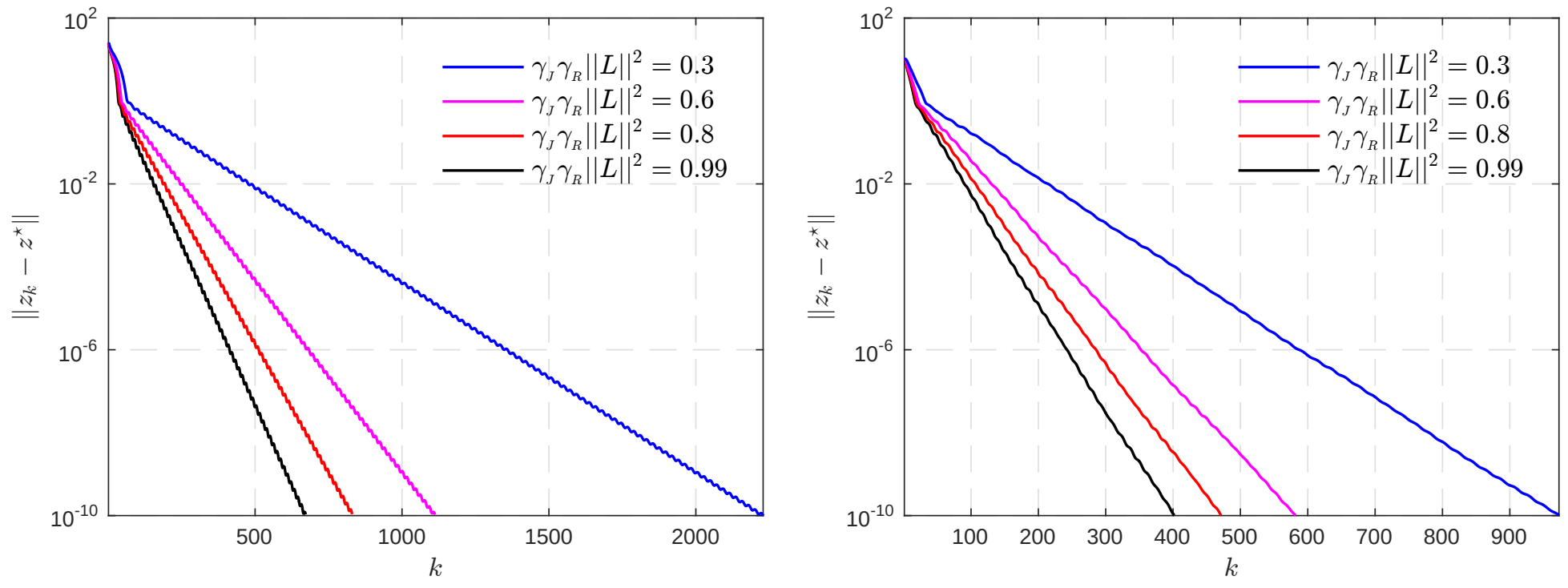

(a) $\ell_1$-norm (b) $\ell_{1,2}$-norm

Figure 3: Comparison of the choice of $\gamma_J, \gamma_R$ when $\theta$ is fixed.

**Fixed $\gamma_J \gamma_R \|L\|^2$** Now we turn to the opposite direction, fix the value of $\gamma_J \gamma_R \|L\|^2$ and then change $\theta$. In the test, we fixed $\gamma_J \gamma_R \|L\|^2 = 0.9$ and $\gamma_J = \gamma_R$, 5 different choices of $\theta$ are considered, which are

$$\theta \in \{0.5, 0.75, 1.0, 2.0\},$$

plus one with Armijo-Goldstein-rule for adaptive update $\theta$. Although there's no convergence guarantee for $\theta = 2.0$, in the tests it converges and we choose to put it here as an illustration of the effects of $\theta$. The result is shown in Figure 4, and we have the following observations

(i) Similar to the previous one, the smaller the value of $\theta$, the slower the iteration converges. Also, the Armijo-Goldstein-rule is the fastest one of all.

(ii) Interestingly, the value of $\theta$ has no impacts to the identification of the iteration.

**Fixed $\theta$ and $\gamma_J \gamma_R$** For the above comparisons, we fix $\gamma_J = \gamma_R$, so for this comparison, we compare the different choices of them. We fix $\theta = 1$ and $\gamma_J \gamma_R \|L\|^2 = 0.99$, then we choose

$$\gamma_R \in \{0.25, 0.5, 1, 2\} \text{ and } \gamma_J = \frac{0.99}{\gamma_R \|L\|^2}.$$

Figure 5 shows the comparison result, we can also have two observations

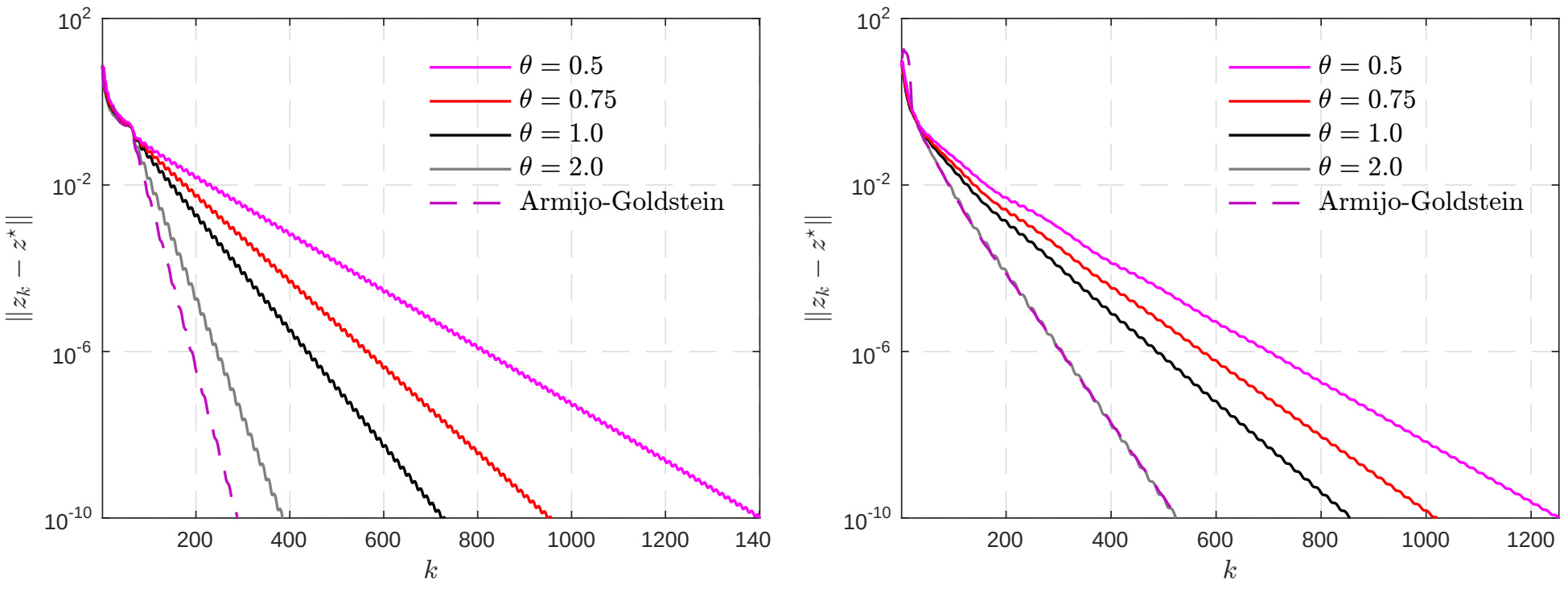

(a) $\ell_1$-norm (b) $\ell_{1,2}$-norm

Figure 4: Comparison of the choice of $\theta$ when $\gamma_J, \gamma_R$ are fixed.

(i) For the $\ell_1$-norm, since both functions are polyhedral, local convergence rate are the same for all choices of $\gamma_R$, see (3.13) for the expression of the rate. The only difference is the identification speed, $\gamma_R = 0.25$ gives the slowest identification, however it uses almost the same number of iterations reaching the given accuracy;

(ii) For the $\ell_{1,2}$-norm, on the other hand, the choice of $\gamma_R$ affects both the identification and local convergence rate. It can be observed that bigger $\gamma_R$ leads to faster local rate, however, it does not mean that the bigger the better. In fact, too big value will slow down the convergence.

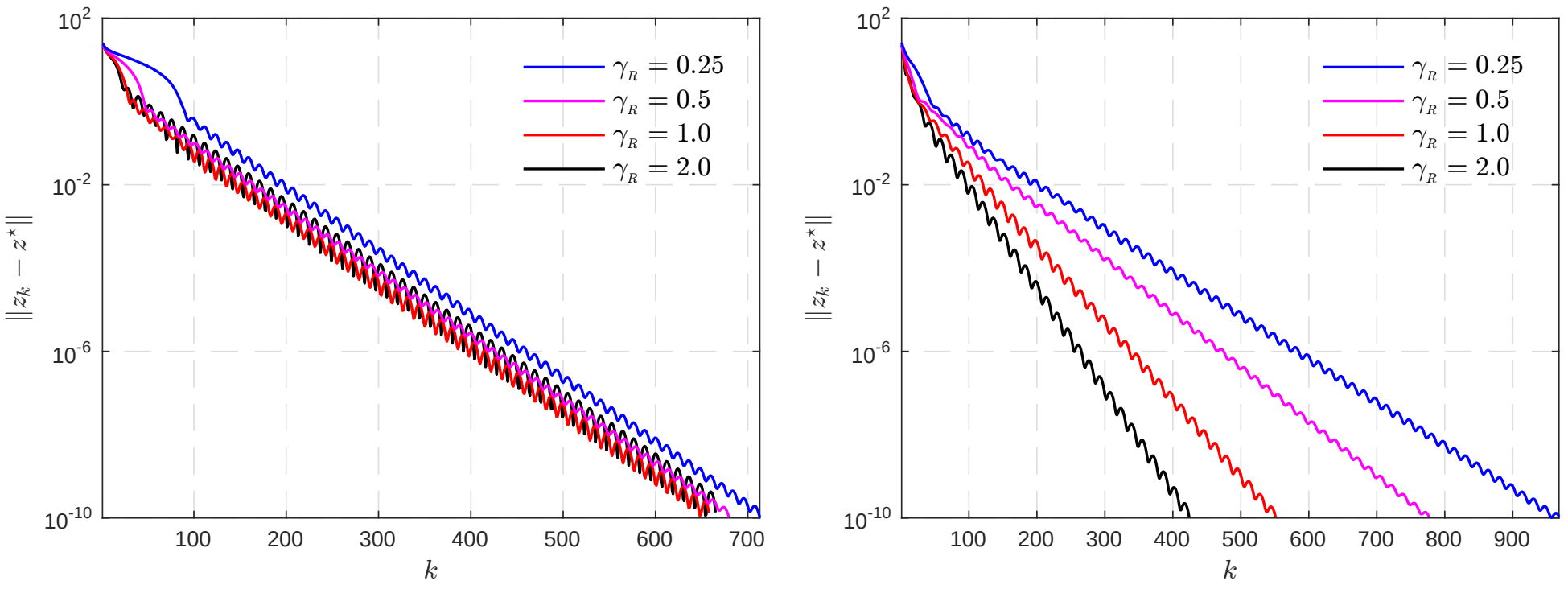

(a) $\ell_1$-norm (b) $\ell_{1,2}$-norm

Figure 5: Comparison of fixed $\theta$ and $\gamma_J\gamma_R$, but varying $\gamma_R$.

To summarize the above comparison, in practice, it is better to chose $\theta$ and $\gamma_J\gamma_R$ as big as possible, moreover, the choices of $\gamma_J$ and $\gamma_R$ should be determined based on the properties of the functions at hand (*i.e.* polyhedral or others).

## 6.5 Oscillation of the method

We dedicate the last part of the numerical experiment to the demonstration of the oscillation behaviour of the Primal–Dual splitting method when dealing with polyhedral functions. As we have seen from the above experiments, oscillation of $\|z_k - z^\star\|$ happens for all examples whose involved functions $R, J^*$ are polyhedral, even for the non-polyhedral $\ell_{1,2}$-norm (for the $\ell_\infty$-norm, due to the fact that the oscillation period is too small compared to the number of iteration, hence it is not visible).

Now to verify our discussion in Section 4, we consider problem (6.2) with $R$ being the $\ell_1$-norm for this illustration, and the result is shown in Figure 6. As revealed in (4.4), the argument of the leading eigenvalue of $M_{\mathrm{PD}} - M_{\mathrm{PD}}^{\infty}$ is controlled by $\theta\gamma_J\gamma_R$, so is the oscillation period. Therefore, the value $\gamma_J\gamma_R$ is tuned such that the oscillation period is an integer, and $\pi/\omega = 12$ for the example we tested. Figure 6 shows graphically the observed oscillation, apparently the oscillation pattern coincides well with the theoretical estimation.

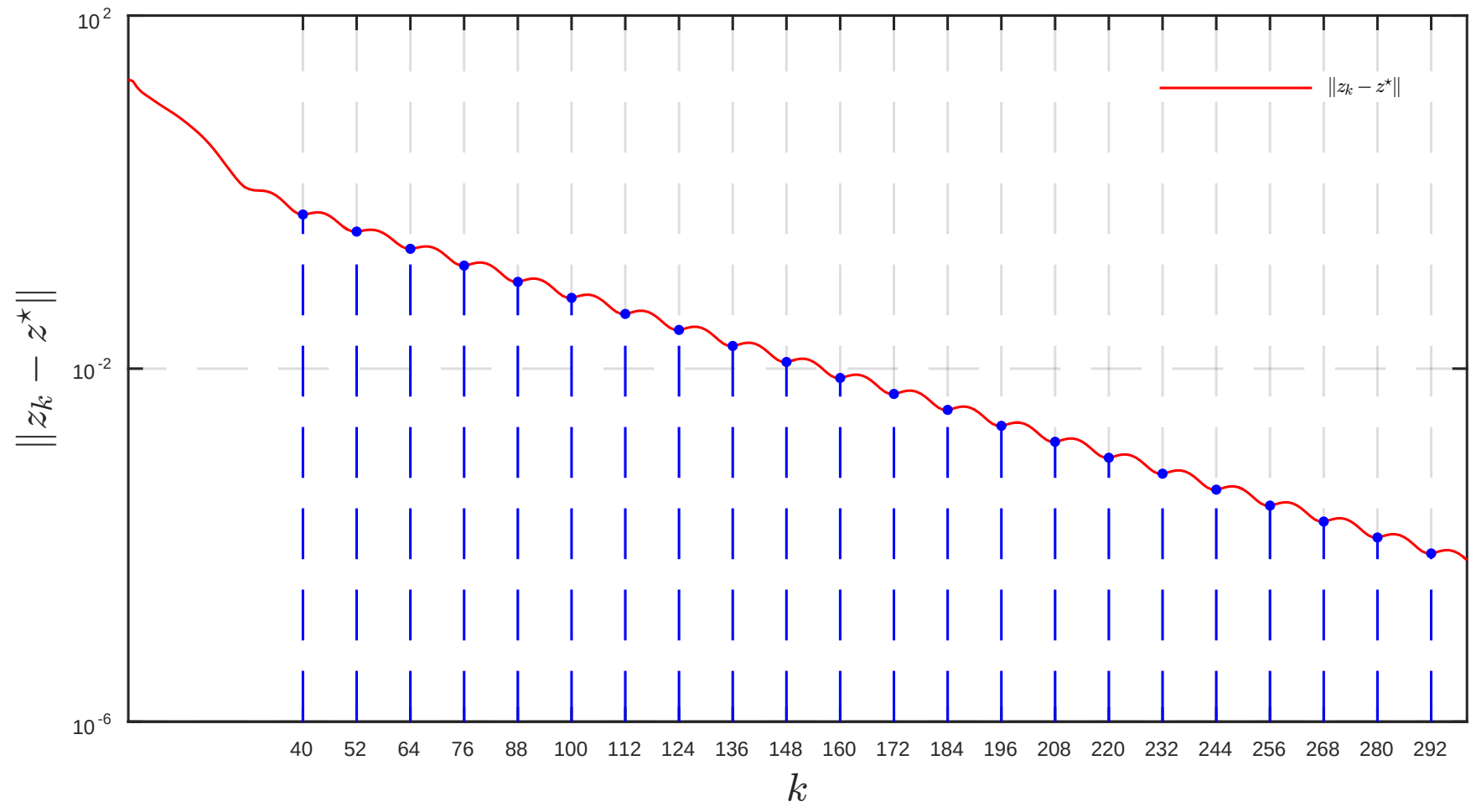

Figure 6: Oscillation behaviour of Primal–Dual splitting method when dealing with polyhedral functions.

# 7 Discussion and conclusion

In this paper, we studied local convergence properties of a class of Primal–Dual splitting methods when the involved functions are moreover partly smooth. In particular, we demonstrated that these methods identify the active manifolds in finite time and then converge locally linearly at a rate that we characterized precisely. We also built connections of the presented result with our previous work on forward–backward splitting method and Douglas–Rachford splitting method/ADMM. Though we focused on one class of Primal–Dual splitting methods, there are other Primal–Dual splitting schemes, such as those in [11, 18, 26], for which our analysis and conclusions can be straightforwardly extended.

# Acknowledgments

This work has been partly supported by the European Research Council (ERC project SIGMA-Vision). JF was partly supported by Institut Universitaire de France. The authors would like to thank Russell Luke for helpful discussions.

# A Proof of Section 3

## A.1 Riemannian Geometry

Let $\mathcal{M}$ be a $C^2$-smooth embedded submanifold of $\mathbb{R}^n$ around a point $x$. With some abuse of terminology, we shall state $C^2$-manifold instead of $C^2$-smooth embedded submanifold of $\mathbb{R}^n$. The natural embedding of a submanifold $\mathcal{M}$ into $\mathbb{R}^n$ permits to define a Riemannian structure and to introduce geodesics on $\mathcal{M}$, and we simply say $\mathcal{M}$ is a Riemannian manifold. We denote respectively $\mathcal{T}_{\mathcal{M}}(x)$ and $\mathcal{N}_{\mathcal{M}}(x)$ the tangent and normal space of $\mathcal{M}$ at point near $x$ in $\mathcal{M}$.

**Exponential map** Geodesics generalize the concept of straight lines in $\mathbb{R}^n$, preserving the zero acceleration characteristic, to manifolds. Roughly speaking, a geodesic is locally the shortest path between two points on $\mathcal{M}$. We denote by $\mathfrak{g}(t;x,h)$ the value at $t \in \mathbb{R}$ of the geodesic starting at $\mathfrak{g}(0;x,h) = x \in \mathcal{M}$ with velocity $\dot{\mathfrak{g}}(t;x,h) = \frac{d\mathfrak{g}}{dt}(t;x,h) = h \in \mathcal{T}_{\mathcal{M}}(x)$ (which is uniquely defined). For every $h \in \mathcal{T}_{\mathcal{M}}(x)$, there exists an interval $I$ around $0$ and a unique geodesic $\mathfrak{g}(t;x,h) : I \to \mathcal{M}$ such that $\mathfrak{g}(0;x,h) = x$ and $\dot{\mathfrak{g}}(0;x,h) = h$. The mapping

$$\mathrm{Exp}_x : \mathcal{T}_{\mathcal{M}}(x) \to \mathcal{M},\ h \mapsto \mathrm{Exp}_x(h) = \mathfrak{g}(1;x,h),$$

is called *Exponential map*. Given $x, x' \in \mathcal{M}$, the direction $h \in \mathcal{T}_{\mathcal{M}}(x)$ we are interested in is such that

$$\mathrm{Exp}_x(h) = x' = \mathfrak{g}(1;x,h).$$

**Parallel translation** Given two points $x, x' \in \mathcal{M}$, let $\mathcal{T}_{\mathcal{M}}(x), \mathcal{T}_{\mathcal{M}}(x')$ be their corresponding tangent spaces. Define

$$\tau : \mathcal{T}_{\mathcal{M}}(x) \to \mathcal{T}_{\mathcal{M}}(x'),$$

the parallel translation along the unique geodesic joining $x$ to $x'$, which is isomorphism and isometry w.r.t. the Riemannian metric.

**Riemannian gradient and Hessian** For a vector $v \in \mathcal{N}_{\mathcal{M}}(x)$, the Weingarten map of $\mathcal{M}$ at $x$ is the operator $\mathfrak{W}_x(\cdot, v) : \mathcal{T}_{\mathcal{M}}(x) \to \mathcal{T}_{\mathcal{M}}(x)$ defined by

$$\mathfrak{W}_x(\cdot, v) = -\mathrm{P}_{\mathcal{T}_{\mathcal{M}}(x)} \mathrm{d}V[h],$$

where $V$ is any local extension of $v$ to a normal vector field on $\mathcal{M}$. The definition is independent of the choice of the extension $V$, and $\mathfrak{W}_x(\cdot, v)$ is a symmetric linear operator which is closely tied to the second fundamental form of $\mathcal{M}$, see [15, Proposition II.2.1].

Let $J$ be a real-valued function which is $C^2$ along the $\mathcal{M}$ around $x$. The covariant gradient of $J$ at $x' \in \mathcal{M}$ is the vector $\nabla_{\mathcal{M}} J(x') \in \mathcal{T}_{\mathcal{M}}(x')$ defined by

$$\langle \nabla_{\mathcal{M}} J(x'),\, h \rangle = \frac{d}{dt} J\big(\mathrm{P}_{\mathcal{M}}(x' + th)\big)\big|_{t=0},\ \forall h \in \mathcal{T}_{\mathcal{M}}(x'),$$

where $\mathrm{P}_{\mathcal{M}}$ is the projection operator onto $\mathcal{M}$. The covariant Hessian of $J$ at $x'$ is the symmetric linear mapping $\nabla^2_{\mathcal{M}} J(x')$ from $\mathcal{T}_{\mathcal{M}}(x')$ to itself which is defined as

$$\langle \nabla^2_{\mathcal{M}} J(x')h,\, h \rangle = \frac{d^2}{dt^2} J\big(\mathrm{P}_{\mathcal{M}}(x' + th)\big)\big|_{t=0},\ \forall h \in \mathcal{T}_{\mathcal{M}}(x'). \tag{A.1}$$

This definition agrees with the usual definition using geodesics or connections [40]. Now assume that $\mathcal{M}$ is a Riemannian embedded submanifold of $\mathbb{R}^n$, and that a function $J$ has a $C^2$-smooth restriction on $\mathcal{M}$. This can be characterized by the existence of a $C^2$-smooth extension (representative) of $J$, *i.e.* a $C^2$-smooth function $\widetilde{J}$ on $\mathbb{R}^n$ such that $\widetilde{J}$ agrees with $J$ on $\mathcal{M}$. Thus, the Riemannian gradient $\nabla_{\mathcal{M}} J(x')$ is also given by

$$\nabla_{\mathcal{M}} J(x') = \mathrm{P}_{\mathcal{T}_{\mathcal{M}}(x')} \nabla \widetilde{J}(x'), \tag{A.2}$$

and $\forall h \in \mathcal{T}_{\mathcal{M}}(x')$, the Riemannian Hessian reads

$$\begin{aligned} \nabla^2_{\mathcal{M}} J(x')h &= \mathrm{P}_{\mathcal{T}_{\mathcal{M}}(x')} \mathrm{d}(\nabla_{\mathcal{M}} J)(x')[h] = \mathrm{P}_{\mathcal{T}_{\mathcal{M}}(x')} \mathrm{d}\big(x' \mapsto \mathrm{P}_{\mathcal{T}_{\mathcal{M}}(x')} \nabla_{\mathcal{M}} \widetilde{J}\big)[h] \\ &= \mathrm{P}_{\mathcal{T}_{\mathcal{M}}(x')} \nabla^2 \widetilde{J}(x')h + \mathfrak{W}_{x'}\big(h, \mathrm{P}_{\mathcal{N}_{\mathcal{M}}(x')} \nabla \widetilde{J}(x')\big), \end{aligned} \tag{A.3}$$

where the last equality comes from [1, Theorem 1]. When $\mathcal{M}$ is an affine or linear subspace of $\mathbb{R}^n$, then obviously $\mathcal{M} = x + \mathcal{T}_{\mathcal{M}}(x)$, and $\mathfrak{W}_{x'}(h, \mathrm{P}_{\mathcal{N}_{\mathcal{M}}(x')} \nabla \widetilde{J}(x')) = 0$, hence (A.3) reduces to

$$\nabla^2_{\mathcal{M}} J(x') = \mathrm{P}_{\mathcal{T}_{\mathcal{M}}(x')} \nabla^2 \widetilde{J}(x') \mathrm{P}_{\mathcal{T}_{\mathcal{M}}(x')}.$$

See [31, 15] for more materials on differential and Riemannian manifolds.

We have the following proposition characterising the parallel translation and the Riemannian Hessian of two close points in $\mathcal{M}$.

**Lemma A.1.** *Let $x, x'$ be two close points in $\mathcal{M}$, denote $\mathcal{T}_{\mathcal{M}}(x), \mathcal{T}_{\mathcal{M}}(x')$ be the tangent spaces of $\mathcal{M}$ at $x, x'$ respectively, and $\tau : \mathcal{T}_{\mathcal{M}}(x') \to \mathcal{T}_{\mathcal{M}}(x)$ be the parallel translation along the unique geodesic joining from $x$ to $x'$, then for the parallel translation we have, given any bounded vector $v \in \mathbb{R}^n$*

$$(\tau \mathrm{P}_{\mathcal{T}_{\mathcal{M}}(x')} - \mathrm{P}_{\mathcal{T}_{\mathcal{M}}(x)})v = o(\|v\|). \tag{A.4}$$

*The Riemannian Taylor expansion of $J \in C^2(\mathcal{M})$ at $x$ for $x'$ reads,*

$$\tau \nabla_{\mathcal{M}} J(x') = \nabla_{\mathcal{M}} J(x) + \nabla^2_{\mathcal{M}} J(x) \mathrm{P}_{\mathcal{T}_{\mathcal{M}}(x)}(x' - x) + o(\|x' - x\|). \tag{A.5}$$

**Proof.** See [36, Lemma B.1 and B.2]. □

**Lemma A.2.** *Let $\mathcal{M}$ be a $C^2$-smooth manifold, $\bar{x} \in \mathcal{M}$, $R \in \mathrm{PSF}_{\bar{x}}(\mathcal{M})$ and $\bar{u} \in \partial R(\bar{x})$. Let $\widetilde{R}$ be a smooth representative of $R$ on $\mathcal{M}$ near $x$, then given any $h \in T_{\bar{x}}$,*

(i) *when $\mathcal{M}$ is a general smooth manifold, if there holds $\bar{u} \in \mathrm{ri}\big(\partial R(\bar{x})\big)$, define the function $\overline{R}(x) = R(x) - \langle x,\, \bar{u}\rangle$, then*

$$\langle \mathrm{P}_{T_{\bar{x}}} \nabla^2_{\mathcal{M}} \overline{R}(\bar{x}) \mathrm{P}_{T_{\bar{x}}} h,\, h\rangle \geq 0. \tag{A.6}$$

(ii) *if $\mathcal{M}$ is affine/linear, then we have directly,*

$$\langle \mathrm{P}_{T_{\bar{x}}} \nabla^2_{\mathcal{M}} \widetilde{R}(\bar{x}) \mathrm{P}_{T_{\bar{x}}} h,\, h\rangle \geq 0. \tag{A.7}$$

**Proof.** See [36, Lemma 4.3]. □

## A.2 Proofs of main theorems

**Proof of Theorem 3.3.**

(i) From the iteration scheme (1.2), for the updating of $x_{k+1}$, we have for $x_{k+1}$,

$$\begin{aligned} & x_{k+1} = \mathrm{prox}_{\gamma_R R}\big(x_k - \gamma_R \nabla F(x_k)) - \gamma_R L^* v_k\big) \\ \Longleftrightarrow\ & \frac{1}{\gamma_R}\big(x_k - \gamma_R \nabla F(x_k) - \gamma_R L^* v_k - x_{k+1}\big) \in \partial R(x_{k+1}), \end{aligned}$$

and then

$$\begin{aligned} \mathrm{dist}\big(-L^* v^\star - \nabla F(x^\star), \partial R(x_{k+1})\big) &\leq \| - L^* v^\star - \nabla F(x^\star) - \tfrac{1}{\gamma_R}(x_k - \gamma_R \nabla F(x_k) - \gamma_R L^* v_k - x_{k+1})\| \\ &\leq \big(\frac{1}{\gamma_R} + \frac{1}{\beta_F}\big)\|x_k - x^\star\| + \|L\|\|v_k - v^\star\| \to 0. \end{aligned}$$

Similarly for $v_{k+1}$

$$\begin{aligned} & v_{k+1} = \mathrm{prox}_{\gamma_J J^*}\big(v_k - \gamma_J \nabla G^*(v_k) + \gamma_J L\bar{x}_{k+1}\big) \\ \Longleftrightarrow\ & \frac{1}{\gamma_J}\big(v_k - \gamma_J \nabla G^*(v_k) + \gamma_J L\bar{x}_{k+1} - v_{k+1}\big) \in \partial J^*(v_{k+1}), \end{aligned}$$

then we have

$$\begin{aligned} & \mathrm{dist}\big(Lx^\star - \nabla G^*(v^\star), \partial J^*(v_{k+1})\big) \\ & \leq \|Lx^\star - \nabla G^*(v^\star) - \tfrac{1}{\gamma_J}(v_k - \gamma_J \nabla G^*(v_k) + \gamma_J L\bar{x}_{k+1} - v_{k+1})\| \\ & \leq \frac{1}{\gamma_J}\|v_k - v_{k+1}\| + \|\nabla G^*(v_k) - \nabla G^*(v^\star)\| + \|L\|\|\bar{x}_{k+1} - x^\star\| \\ & \leq \big(\frac{1}{\gamma_J} + \frac{1}{\beta_G}\big)\|v_k - v^\star\| + \|L\|\big((1+\theta)\|x_{k+1} - x^\star\| + \theta\|x_k - x^\star\|\big) \to 0. \end{aligned}$$

By assumption, $J^* \in \Gamma_0(\mathbb{R}^m)$, $R \in \Gamma_0(\mathbb{R}^n)$, hence they are sub-differentially continuous at every point in their respective domains [44, Example 13.30], and in particular at $v^\star$ and $x^\star$. It then follows that $J^*(v_k) \to J^*(v^\star)$ and $R(x_k) \to R(x^\star)$. Altogether with the non-degeneracy condition (ND), shows that the conditions of [28, Theorem 5.3] are fulfilled for $\langle -Lx^\star + \nabla G^*(v^\star),\, \cdot\rangle + J^*$ and $\langle L^* v^\star + \nabla F(x^\star),\, \cdot\rangle + R$, and the finite identification claim follows.

(ii) (a) In this case, $\mathcal{M}^R_{x^\star}$ is an affine subspace, *i.e.* $\mathcal{M}^R_{x^\star} = x^\star + T^R_{x^\star}$, it is straight to have $\bar{x}_k \in \mathcal{M}^R_{x^\star}$. Then since $R$ is partly smooth at $x^\star$ relative to $\mathcal{M}^R_{x^\star}$, the sharpness property holds at all nearby points in $\mathcal{M}^R_{x^\star}$ [33, Proposition 2.10]. Thus for $k$ large enough, i.e. $x_k$ sufficiently close to $x^\star$ on $\mathcal{M}^R_{x^\star}$, we have indeed $\mathcal{T}_{x_k}(\mathcal{M}^R_{x^\star}) = T^R_{x^\star} = T^R_{x_k}$ as claimed.

(b) Similar to (ii)(a).

(c) It is immediate to verify that a locally polyhedral function around $x^\star$ is indeed partly smooth relative to the affine subspace $x^\star + T^R_{x^\star}$, and thus, the first claim follows from (ii)(a). For the rest, it is sufficient to observe that by polyhedrality, for any $x \in \mathcal{M}^R_{x^\star}$ near $x^\star$, $\partial R(x) = \partial R(x^\star)$. Therefore, combining local normal sharpness [33, Proposition 2.10] and Lemma A.2 yields the second conclusion.

(d) Similar to (ii)(c). □

**Proof of Proposition 3.6.** From the update of $x_k$ in (1.2), we have

$$x_k - \gamma_R \nabla F(x_k) - \gamma_R L^* v_k - x_{k+1} \in \gamma_R \partial R(x_{k+1}),$$
$$-\gamma_R \nabla F(x^\star) - \gamma_R L^* v^\star \in \gamma_R \partial R(x^\star).$$

Denote $\tau^R_k$ the parallel translation from $T^R_{x_k}$ to $T^R_{x^\star}$. Then project on to corresponding tangent spaces and apply parallel translation,

$$\begin{aligned}\gamma_R \tau^R_k \nabla_{\mathcal{M}^R_{x^\star}} R(x_{k+1}) &= \tau^R_k \mathrm{P}_{T^R_{x^\star} x_{k+1}}(x_k - \gamma_R \nabla F(x_k) - \gamma_R L^* v_k - x_{k+1}) \\ &= \mathrm{P}_{T^R_{x^\star}}(x_k - \gamma_R \nabla F(x_k) - \gamma_R L^* v_k - x_{k+1}) \\ &\quad + (\tau^R_k \mathrm{P}_{T^R_{x^\star} x_{k+1}} - \mathrm{P}_{T^R_{x^\star}})(x_k - \gamma_R \nabla F(x_k) - \gamma_R L^* v_k - x_{k+1}), \\ \gamma_R \nabla_{\mathcal{M}^R_{x^\star}} R(x^\star) &= \mathrm{P}_{T^R_{x^\star}}(-\gamma_R \nabla F(x^\star) - \gamma_R L^* v^\star),\end{aligned}$$

which leads to

$$\begin{aligned}&\gamma_R \tau^R_k \nabla_{\mathcal{M}^R_{x^\star}} R(x_{k+1}) - \gamma_R \nabla_{\mathcal{M}^R_{x^\star}} R(x^\star) \\ &= \mathrm{P}_{T^R_{x^\star}}\big((x_k - \gamma_R \nabla F(x_k) - \gamma_R L^* v_k - x_{k+1}) - (x^\star - \gamma_R \nabla F(x^\star) - \gamma_R L^* v^\star - x^\star)\big) \\ &\quad + \underbrace{(\tau^R_k \mathrm{P}_{T^R_{x^\star} x_{k+1}} - \mathrm{P}_{T^R_{x^\star}})(-\gamma_R \nabla F(x^\star) - \gamma_R L^* v^\star)}_{\textbf{Term 1}} \\ &\quad + \underbrace{(\tau^R_k \mathrm{P}_{T^R_{x^\star} x_{k+1}} - \mathrm{P}_{T^R_{x^\star}})\big((x_k - \gamma_R \nabla F(x_k) - \gamma_R L^* v_k - x_{k+1}) + (\gamma_R \nabla F(x^\star) + \gamma_R L^* v^\star)\big)}_{\textbf{Term 2}}.\end{aligned} \tag{A.8}$$

Moving **Term 1** to the other side leads to

$$\begin{aligned}&\gamma_R \tau^R_k \nabla_{\mathcal{M}^R_{x^\star}} R(x_{k+1}) - \gamma_R \nabla_{\mathcal{M}^R_{x^\star}} R(x^\star) - (\tau^R_k \mathrm{P}_{T^R_{x^\star} x_{k+1}} - \mathrm{P}_{T^R_{x^\star}})(-\gamma_R \nabla F(x^\star) - \gamma_R L^* v^\star) \\ &= \gamma_R \tau^R_k \big(\nabla_{\mathcal{M}^R_{x^\star}} R(x_{k+1}) + (L^* v^\star + \nabla F(x^\star))\big) - \gamma_R \big(\nabla_{\mathcal{M}^R_{x^\star}} R(x^\star) + (L^* v^\star + \nabla F(x^\star))\big) \\ &= \gamma_R \mathrm{P}_{T^R_{x^\star}} \nabla^2_{\mathcal{M}^R_{x^\star}} \overline{R}(x^\star) \mathrm{P}_{T^R_{x^\star}}(x_{k+1} - x^\star) + o(\|x_{k+1} - x^\star\|),\end{aligned}$$

where Lemma A.2 is applied. Since $x_{k+1} = \mathrm{prox}_{\gamma_R R}(x_k - \gamma_R \nabla F(x_k) - \gamma_R L^* v_k)$, $\mathrm{prox}_{\gamma_R R}$ is firmly non-expansive and $\mathrm{Id}_n - \gamma_R \nabla F$ is non-expansive (under the parameter setting), then

$$\begin{aligned}&\|(x_k - \gamma_R \nabla F(x_k) - \gamma_R L^* v_k - x_{k+1}) - (x^\star - \gamma_R \nabla F(x^\star) - \gamma_R L^* v^\star - x^\star)\| \\ &\le \|(\mathrm{Id}_n - \gamma_R \nabla F)(x_k) - (\mathrm{Id}_n - \gamma_R \nabla F)(x^\star)\| + \gamma_R \|L^* v_k - L^* v^\star\| \\ &\le \|x_k - x^\star\| + \gamma_R \|L\| \|v_k - v^\star\|.\end{aligned} \tag{A.9}$$

Therefore, for **Term 2**, owing to Lemma A.1, we have

$$\begin{aligned}&(\tau^R_k \mathrm{P}_{T^R_{x^\star} x_{k+1}} - \mathrm{P}_{T^R_{x^\star}})\big((x_k - \gamma_R \nabla F(x_k) - \gamma_R L^* v_k - x_{k+1}) - (x^\star - \gamma_R \nabla F(x^\star) - \gamma_R L^* v^\star - x^\star)\big) \\ &= o(\|x_k - x^\star\| + \gamma_R \|L\| \|v_k - v^\star\|).\end{aligned}$$

Therefore, from (A.8), and apply $x_k - x^\star = \mathrm{P}_{T^R_{x^\star}}(x_k - x^\star) + o(x_k - x^\star)$ [34, Lemma 5.1] to $(x_{k+1} - x^\star)$ and $(x_k - x^\star)$, we get

$$\begin{aligned}&(\mathrm{Id}_n + H_{\overline{R}})(x_{k+1} - x^\star) + o(\|x_{k+1} - x^\star\|) \\ &= (x_k - x^\star) - \gamma_R \mathrm{P}_{T^R_{x^\star}}(\nabla F(x_k) - \nabla F(x^\star)) - \gamma_R \mathrm{P}_{T^R_{x^\star}} L^*(v_k - v^\star) + o(\|x_k - x^\star\| + \gamma_R \|L\| \|v_k - v^\star\|).\end{aligned}$$

Then apply Taylor expansion to $\nabla F$, and apply [34, Lemma 5.1] to $(v_k - v^\star)$,

$$\begin{aligned}&(\mathrm{Id}_n + H_{\overline{R}})(x_{k+1} - x^\star)\\&= (\mathrm{Id}_n - \gamma_R H_F)(x_k - x^\star) - \gamma_R \overline{L}^*(v_k - v^\star) + o(\|x_k - x^\star\| + \gamma_R\|L\|\|v_k - v^\star\|).\end{aligned}\tag{A.10}$$

Then invert $(\mathrm{Id}_n + H_{\overline{R}})$ and apply [34, Lemma 5.1], we get

$$x_{k+1} - x^\star = W_{\overline{R}}\overline{H}_F(x_k - x^\star) - \gamma_R W_{\overline{R}}\overline{L}^*(v_k - v^\star) + o(\|x_k - x^\star\| + \gamma_R\|L\|\|v_k - v^\star\|).\tag{A.11}$$

Now from the update of $v_{k+1}$

$$\begin{aligned}v_k - \gamma_J \nabla G^*(v_k) + \gamma_J L\bar{x}_{k+1} - v_{k+1} &\in \gamma_J \partial J^*(v_{k+1}),\\ v^\star - \gamma_J \nabla G^*(v^\star) + \gamma_J L x^\star - v^\star &\in \gamma_J \partial J^*(v^\star).\end{aligned}$$

Denote $\tau^{J^*}_{k+1}$ the parallel translation from $T^{J^*}_{v_{k+1}}$ to $T^{J^*}_{v^\star}$, then

$$\begin{aligned}\gamma_J \tau^{J^*}_{k+1}\nabla_{\mathcal{M}^{J^*}_{v^\star}} J^*(v_{k+1}) &= \tau^{J^*}_{k+1}\mathrm{P}_{T^{J^*}_{v_{k+1}}}(v_k - \gamma_J\nabla G^*(v_k) + \gamma_J L\bar{x}_{k+1} - v_{k+1})\\&= \mathrm{P}_{T^{J^*}_{v^\star}}(v_k - \gamma_J\nabla G^*(v_k) + \gamma_J L\bar{x}_{k+1} - v_{k+1})\\&\quad + (\tau^{J^*}_{k+1}\mathrm{P}_{T^{J^*}_{v_{k+1}}} - \mathrm{P}_{T^{J^*}_{v^\star}})(v_k - \gamma_J\nabla G^*(v_k) + \gamma_J L\bar{x}_{k+1} - v_{k+1}),\\ \gamma_J\nabla_{\mathcal{M}^{J^*}_{v^\star}} J^*(v^\star) &= \mathrm{P}_{T^{J^*}_{v^\star}}(v^\star - \gamma_J\nabla G^*(v^\star) + \gamma_J L x^\star - v^\star)\end{aligned}$$

which leads to

$$\begin{aligned}&\gamma_J \tau^{J^*}_{k+1}\nabla_{\mathcal{M}^{J^*}_{v^\star}} J^*(v_{k+1}) - \gamma_J\nabla_{\mathcal{M}^{J^*}_{v^\star}} J^*(v^\star)\\&= \mathrm{P}_{T^{J^*}_{v^\star}}\big((v_k - \gamma_J\nabla G^*(v_k) + \gamma_J L\bar{x}_{k+1} - v_{k+1}) - (v^\star - \gamma_J\nabla G^*(v^\star) + \gamma_J Lx^\star - v^\star)\big)\\&\quad + (\tau^{J^*}_{k+1}\mathrm{P}_{T^{J^*}_{v_{k+1}}} - \mathrm{P}_{T^{J^*}_{v^\star}})(v_k - \gamma_J\nabla G^*(v_k) + \gamma_J L\bar{x}_{k+1} - v_{k+1})\\&= \mathrm{P}_{T^{J^*}_{v^\star}}\big((v_k - \gamma_J\nabla G^*(v_k) + \gamma_J L\bar{x}_{k+1} - v_{k+1}) - (v^\star - \gamma_J\nabla G^*(v^\star) + \gamma_J Lx^\star - v^\star)\big)\\&\quad + \underbrace{(\tau^{J^*}_{k+1}\mathrm{P}_{T^{J^*}_{v_{k+1}}} - \mathrm{P}_{T^{J^*}_{v^\star}})(\gamma_J Lx^\star - \gamma_J\nabla G^*(v^\star))}_{\textbf{Term 3}}\\&\quad + \underbrace{(\tau^{J^*}_{k+1}\mathrm{P}_{T^{J^*}_{v_{k+1}}} - \mathrm{P}_{T^{J^*}_{v^\star}})\big((v_k - \gamma_J\nabla G^*(v_k) + \gamma_J L\bar{x}_{k+1} - v_{k+1}) + \gamma_J(\nabla G^*(v^\star) - Lx^\star)\big)}_{\textbf{Term 4}}.\end{aligned}\tag{A.12}$$

Similarly to the previous analysis, for **Term 3**, move to the lefthand side of the inequality and apply Lemma A.2,

$$\begin{aligned}&\gamma_J \tau^{J^*}_{k+1}\nabla_{\mathcal{M}^{J^*}_{v^\star}} J^*(v_{k+1}) - \gamma_J\nabla_{\mathcal{M}^{J^*}_{v^\star}} J^*(v^\star) - (\tau^{J^*}_{k+1}\mathrm{P}_{T^{J^*}_{v_{k+1}}} - \mathrm{P}_{T^{J^*}_{v^\star}})(\gamma_J Lx^\star - \gamma_J\nabla G^*(v^\star))\\&= \gamma_J\tau^{J^*}_{k+1}\big(\nabla_{\mathcal{M}^{J^*}_{v^\star}} J^*(v_{k+1}) - (Lx^\star - \nabla G^*(v^\star))\big) - \gamma_J\big(\nabla_{\mathcal{M}^{J^*}_{v^\star}} J^*(v^\star) - (Lx^\star - \nabla G^*(v^\star))\big)\\&= \gamma_J \mathrm{P}_{T^{J^*}_{v^\star}}\nabla^2_{\mathcal{M}^{J^*}_{v^\star}}\overline{J}^*(v^\star)\mathrm{P}_{T^{J^*}_{v^\star}}(v_{k+1} - v^\star) + o(\|v_{k+1} - v^\star\|).\end{aligned}$$

Since $\theta \le 1$, we have

$$\begin{aligned}\|\bar{x}_{k+1} - x^\star\| &\le (1+\theta)\|x_{k+1} - x^\star\| + \theta\|x_k - x^\star\|\\&\le 2(\|x_k - x^\star\| + \gamma_R\|L\|\|v_k - v^\star\|) + \|x_k - x^\star\|\\&= 3\|x_k - x^\star\| + 2\gamma_R\|L\|\|v_k - v^\star\|.\end{aligned}$$

Then for **Term 4**, since $\gamma_J\gamma_R\|L^2\| < 1$, $\mathrm{prox}_{\gamma_J J^*}$ is firmly non-expansive and $\mathrm{Id}_m - \gamma_J\nabla G^*$ is non-expansive, we have

$$\begin{aligned}&(\tau^{J^*}_{k+1}\mathrm{P}_{T^{J^*}_{v_{k+1}}} - \mathrm{P}_{T^{J^*}_{v^\star}})\big((v_k - \gamma_J\nabla G^*(v_k) + \gamma_J L\bar{x}_{k+1} - v_{k+1}) - (v^\star - \gamma_J\nabla G^*(v^\star) + \gamma_J Lx^\star - v^\star)\big)\\&= o(\|v_k - v^\star\| + \gamma_J\|L\|\|x_k - x^\star\|).\end{aligned}$$

Therefore, from (A.12), apply [34, Lemma 5.1] to $(v_{k+1} - v^\star)$ and $(v_k - v^\star)$, we get

$$(\mathrm{Id}_m + H_{\overline{J^*}})(v_{k+1} - v^\star) = (\mathrm{Id}_m - \gamma_J H_{G^*})(v_k - v^\star) + \gamma_J \overline{L}(\bar{x}_{k+1} - x^\star) + o(\|v_k - v^\star\| + \gamma_J \|L\|\|x_k - x^\star\|). \tag{A.13}$$

Then similar to (A.11), we get from (A.13)

$$\begin{aligned}
v_{k+1} - v^\star &= W_{\overline{J^*}} \overline{H}_{G^*}(v_k - v^\star) + \gamma_J W_{\overline{J^*}} \overline{L}(\bar{x}_{k+1} - x^\star) + o(\|v_k - v^\star\| + \gamma_J \|L\|\|x_k - x^\star\|) \\
&= W_{\overline{J^*}} \overline{H}_{G^*}(v_k - v^\star) + (1+\theta)\gamma_J W_{\overline{J^*}} \overline{L}(x_{k+1} - x^\star) - \theta\gamma_J W_{\overline{J^*}} \overline{L}(x_k - x^\star) \\
&\quad + o(\|v_k - v^\star\| + \gamma_J \|L\|\|x_k - x^\star\|) \\
&= W_{\overline{J^*}} \overline{H}_{G^*}(v_k - v^\star) - \theta\gamma_J W_{\overline{J^*}} \overline{L}(x_k - x^\star) \\
&\quad + (1+\theta)\gamma_J W_{\overline{J^*}} \overline{L}\big(W_{\overline{R}} \overline{H}_F (x_k - x^\star) - \gamma_R W_{\overline{R}} \overline{L}^*(v_k - v^\star)\big) \\
&\quad + o(\|x_k - x^\star\| + \gamma_R \|L\|\|v_k - v^\star\|) + o(\|v_k - v^\star\| + \gamma_J \|L\|\|x_k - x^\star\|) \\
&= \big(W_{\overline{J^*}} \overline{H}_{G^*} - (1+\theta)\gamma_J \gamma_R W_{\overline{J^*}} \overline{L} W_{\overline{R}} \overline{L}^*\big)(v_k - v^\star) \\
&\quad + \big((1+\theta)\gamma_J W_{\overline{J^*}} \overline{L} W_{\overline{R}} \overline{H}_F - \theta\gamma_J W_{\overline{J^*}} \overline{L}\big)(x_k - x^\star) \\
&\quad + o(\|x_k - x^\star\| + \gamma_R \|L\|\|v_k - v^\star\|) + o(\|v_k - v^\star\| + \gamma_J \|L\|\|x_k - x^\star\|).
\end{aligned} \tag{A.14}$$

Now we consider the small $o$-terms. For the 2 small $o$-terms in (A.10) and (A.13). First, let $a_1, a_2$ be two constants, then we have

$$|a_1| + |a_2| = \sqrt{(|a_1| + |a_2|)^2} \leq \sqrt{2(a_1^2 + a_2^2)} = \sqrt{2}\left\| \begin{pmatrix} a_1 \\ a_2 \end{pmatrix} \right\|.$$

Denote $b = \max\{1, \gamma_J\|L\|, \gamma_R\|L\|\}$, then

$$\begin{aligned}
&(\|v_k - v^\star\| + \gamma_J \|L\|\|x_k - x^\star\|) + (\|x_k - x^\star\| + \gamma_R \|L\|\|v_k - v^\star\|) \\
&\leq 2b(\|x_k - x^\star\| + \|v_k - v^\star\|) \leq 2\sqrt{2}b \left\| \begin{pmatrix} x_k - x^\star \\ v_k - v^\star \end{pmatrix} \right\|.
\end{aligned}$$

Combining this with (A.11) and (A.14), and ignoring the constants of the small $o$-term leads to the claimed result. □

**Proof of Proposition 3.8.**

(i) When $\theta = 1$, $M_{\mathrm{PD}}$ becomes

$$M_{\mathrm{PD}} = \begin{bmatrix} W_{\overline{R}} \overline{H}_F & -\gamma_R W_{\overline{R}} \overline{L}^* \\ 2\gamma_J W_{\overline{J^*}} \overline{L} W_{\overline{R}} \overline{H}_F - \gamma_J W_{\overline{J^*}} \overline{L} & W_{\overline{J^*}} \overline{H}_{G^*} - 2\gamma_R \gamma_J W_{\overline{J^*}} \overline{L} W_{\overline{R}} \overline{L}^* \end{bmatrix} \tag{A.15}$$

Next we show that $M_{\mathrm{PD}}$ is *averaged non-expansive*.

First define the following matrices

$$\boldsymbol{A} = \begin{bmatrix} H_R/\gamma_R & \overline{L}^* \\ -\overline{L} & H_{J^*}/\gamma_J \end{bmatrix}, \quad \boldsymbol{B} = \begin{bmatrix} H_F & 0 \\ 0 & H_{G^*} \end{bmatrix}, \quad \boldsymbol{\mathcal{V}} = \begin{bmatrix} \mathrm{Id}_n/\gamma_R & -\overline{L}^* \\ -\overline{L} & \mathrm{Id}_m/\gamma_J \end{bmatrix}, \tag{A.16}$$

where we have $\boldsymbol{A}$ is maximal monotone [11], $\boldsymbol{B}$ is $\min\{\beta_F, \beta_G\}$-cocoercive, and $\boldsymbol{\mathcal{V}}$ is $\nu$-positive definite with $\nu = (1 - \sqrt{\gamma_J \gamma_R \|L\|^2}) \min\{\frac{1}{\gamma_R}, \frac{1}{\gamma_J}\}$.

Now we have

$$\boldsymbol{\mathcal{V}} + \boldsymbol{A} = \begin{bmatrix} \frac{\mathrm{Id}_n + H_R}{\gamma_R} & 0 \\ -2\overline{L} & \frac{\mathrm{Id}_m + H_{J^*}}{\gamma_J} \end{bmatrix} \Longrightarrow (\boldsymbol{\mathcal{V}} + \boldsymbol{A})^{-1} = \begin{bmatrix} \gamma_R W_{\overline{R}} & 0 \\ 2\gamma_J \gamma_R W_{\overline{J^*}} \overline{L} W_{\overline{R}} & \gamma_J W_{\overline{J^*}} \end{bmatrix},$$

and

$$\boldsymbol{\mathcal{V}} - \boldsymbol{B} = \begin{bmatrix} \frac{1}{\gamma_R} \overline{H}_F & -\overline{L}^* \\ -\overline{L} & \frac{1}{\gamma_J} \overline{H}_{G^*} \end{bmatrix}.$$

As a result, we get

$$(\mathcal{V}+\boldsymbol{A})^{-1}(\mathcal{V}-\boldsymbol{B}) = \begin{bmatrix} \gamma_R W_{\overline{R}} & 0 \\ 2\gamma_J\gamma_R W_{\overline{J^*}}\overline{L}W_{\overline{R}} & \gamma_J W_{\overline{J^*}} \end{bmatrix} \begin{bmatrix} \frac{1}{\gamma_R}\overline{H}_F & -\overline{L}^* \\ -\overline{L} & \frac{1}{\gamma_J}\overline{H}_{G^*} \end{bmatrix}$$
$$= \begin{bmatrix} W_{\overline{R}}\overline{H}_F & -\gamma_R W_{\overline{R}}\overline{L}^* \\ 2\gamma_J W_{\overline{J^*}}\overline{L}W_{\overline{R}}\overline{H}_F - \gamma_J W_{\overline{J^*}}\overline{L} & W_{\overline{J^*}}\overline{H}_{G^*} - 2\gamma_J\gamma_R W_{\overline{J^*}}\overline{L}W_{\overline{R}}\overline{L}^* \end{bmatrix},$$

which is exactly (A.15).

From Lemma 2.5 we know that $M_{\text{PD}} : \mathcal{K}_{\mathcal{V}} \to \mathcal{K}_{\mathcal{V}}$ is *averaged non-expansive*, hence it is convergent [6]. Then we have the induced matrix norm

$$\lim_{k\to\infty} \|M_{\text{PD}}^k - M_{\text{PD}}^\infty\|_{\mathcal{V}} = \lim_{k\to\infty} \|M_{\text{PD}} - M_{\text{PD}}^\infty\|_{\mathcal{V}}^k = 0.$$

Since we are in the finite dimensional space and $\mathcal{V}$ is an isomorphism, then the above limit implies that

$$\lim_{k\to\infty} \|M_{\text{PD}} - M_{\text{PD}}^\infty\|^k = 0,$$

which means that $\rho(M_{\text{PD}} - M_{\text{PD}}^\infty) < 1$. The rest of the proof is classical using the spectral radius formula, see *e.g.* [5, Theorem 2.12(i)].

(ii) When $R$ and $J^*$ are locally polyhedral, then $W_{\overline{R}} = \text{Id}_n, W_{\overline{J^*}} = \text{Id}_m$, altogether with $F = 0, G = 0$, for any $\theta \in [0,1]$, we have

$$M_{\text{PD}} = \begin{bmatrix} \text{Id}_n & -\gamma_R\overline{L}^* \\ \gamma_J\overline{L} & \text{Id}_m - \gamma_R\gamma_J(1+\theta)\overline{L}\overline{L}^* \end{bmatrix}. \tag{A.17}$$

With the SVD of $\overline{L}$, for $M_{\text{PD}}$, we have

$$\begin{aligned} M_{\text{PD}} &= \begin{bmatrix} \text{Id}_n & -\gamma_R\overline{L}^* \\ \gamma_J\overline{L} & \text{Id}_m - (1+\theta)\gamma_R\gamma_J\overline{L}\overline{L}^* \end{bmatrix} \\ &= \begin{bmatrix} YY^* & -\gamma_R Y\Sigma_{\overline{L}}^* X^* \\ \gamma_J X\Sigma_{\overline{L}}Y^* & XX^* - (1+\theta)\gamma_R\gamma_J X\Sigma_{\overline{L}}^2 X^* \end{bmatrix} \\ &= \begin{bmatrix} Y & \\ & X \end{bmatrix} \underbrace{\begin{bmatrix} \text{Id}_n & -\gamma_R\Sigma_L^* \\ \gamma_J\Sigma_L & \text{Id}_m - (1+\theta)\gamma_R\gamma_J\Sigma_L^2 \end{bmatrix}}_{M_\Sigma} \begin{bmatrix} Y^* & \\ & X^* \end{bmatrix}. \end{aligned} \tag{A.18}$$

Since we assume that $\text{rank}(\overline{L}) = l \leq p$, then $\Sigma_L$ can be represented as

$$\Sigma_L = \begin{bmatrix} \Sigma_l & 0_{l,n-l} \\ 0_{m-l,l} & 0_{m-l,n-l} \end{bmatrix},$$

where $\Sigma_l = (\sigma_j)_{j=1,\ldots,l}$. Back to $M_\Sigma$, we have

$$M_\Sigma = \begin{bmatrix} \text{Id}_l & 0_{l,n-l} & -\gamma_R\Sigma_l & 0_{l,m-l} \\ 0_{n-l,l} & \text{Id}_{n-l} & 0_{n-l,l} & 0_{n-l,m-l} \\ \gamma_J\Sigma_l & 0_{l,n-l} & \text{Id}_l - (1+\theta)\gamma_R\gamma_J\Sigma_l^2 & 0_{l,m-l} \\ 0_{m-l,l} & 0_{m-l,n-l} & 0_{m-l,l} & \text{Id}_{m-l} \end{bmatrix}.$$

Let's study the eigenvalues of $M_\Sigma$,

$$\begin{aligned} |M_\Sigma - \rho\text{Id}_{m+n}| &= \left| \begin{bmatrix} (1-\rho)\text{Id}_l & 0_{l,n-l} & -\gamma_R\Sigma_l & 0_{l,m-l} \\ 0_{n-l,l} & (1-\rho)\text{Id}_{n-l} & 0_{n-l,l} & 0_{n-l,m-l} \\ \gamma_J\Sigma_l & 0_{l,n-l} & (1-\rho)\text{Id}_l - (1+\theta)\gamma_R\gamma_J\Sigma_l^2 & 0_{l,m-l} \\ 0_{m-l,l} & 0_{m-l,n-l} & 0_{m-l,l} & (1-\rho)\text{Id}_{m-l} \end{bmatrix} \right| \\ &= (1-\rho)^{m+n-2l} \left| \begin{bmatrix} (1-\rho)\text{Id}_l & -\gamma_R\Sigma_l \\ \gamma_J\Sigma_l & (1-\rho)\text{Id}_l - (1+\theta)\gamma_R\gamma_J\Sigma_l^2 \end{bmatrix} \right|. \end{aligned}$$

Since $(-\gamma_R \Sigma_l)((1-\rho)\mathrm{Id}_l) = ((1-\rho)\mathrm{Id}_l)(-\gamma_R\Sigma_l)$, then by [45, Theorem 3], we have

$$\begin{aligned}
|M_\Sigma - \rho\mathrm{Id}_{a+b}| &= (1-\rho)^{m+n-2l}\left|\begin{bmatrix}(1-\rho)\mathrm{Id}_l & -\gamma_R\Sigma_l \\ \gamma_J\Sigma_l & (1-\rho)\mathrm{Id}_l-(1+\theta)\gamma_R\gamma_J\Sigma_l^2\end{bmatrix}\right| \\
&= (1-\rho)^{m+n-2l}\left|\left[(1-\rho)\big((1-\rho)\mathrm{Id}_l-(1+\theta)\gamma_R\gamma_J\Sigma_l^2\big)+\gamma_R\gamma_J\Sigma_l\Sigma_l\right]\right| \\
&= (1-\rho)^{m+n-2l}\left|\left[(1-\rho)^2\mathrm{Id}_l-(1-\rho)(1+\theta)\gamma_R\gamma_J\Sigma_l^2+\gamma_R\gamma_J\Sigma_l\Sigma_l\right]\right| \\
&= (1-\rho)^{m+n-2l}\prod\nolimits_{j=1}^{l}\big(\rho^2-(2-(1+\theta)\gamma_J\gamma_R\sigma_j^2)\rho+(1-\theta\gamma_J\gamma_R\sigma_j^2)\big).
\end{aligned}$$

For the eigenvalues $\rho$, clearly, except the 1's, we have for $j=1,...,l$

$$\rho_j = \frac{\big(2-(1+\theta)\gamma_J\gamma_R\sigma_j^2\big) \pm \sqrt{(1+\theta)^2\gamma_J^2\gamma_R^2\sigma_j^4-4\gamma_J\gamma_R\sigma_j^2}}{2}.$$

Since $\gamma_J\gamma_R\sigma_j^2 \le \gamma_J\gamma_R\|L\|^2 < 1$, then $\rho_j$ are complex and

$$|\rho_j| = \frac{1}{2}\sqrt{\big(2-(1+\theta)\gamma_J\gamma_R\sigma_j^2\big)^2-\big((1+\theta)^2\gamma_J^2\gamma_R^2\sigma_j^4-4\gamma_J\gamma_R\sigma_j^2\big)} = \sqrt{1-\theta\gamma_J\gamma_R\sigma_j^2} < 1.$$

As a result, we also obtain the $M_{\mathrm{PD}}^\infty$, which reads

$$M_{\mathrm{PD}}^\infty = \begin{bmatrix} Y & \\ & X\end{bmatrix}\begin{bmatrix}0_l & & & \\ & \mathrm{Id}_{n-l} & & \\ & & 0_l & \\ & & & \mathrm{Id}_{m-l}\end{bmatrix}\begin{bmatrix}Y^* & \\ & X^*\end{bmatrix}.$$

If $n=m$ and moreover $L=\mathrm{Id}_n$, then $(\sigma_j)_{j=1,\dots,p}$ corresponds to *cosine* the *principal angles* between the tangent spaces $T_{x^\star}^R$ and $T_{v^\star}^{J^*}$ [24]. □

**Proof of Corollary 3.10.**

(i) From the local fixed-point iteration (3.10), we have

$$z_{k+1}-z^\star = M_{\mathrm{PD}}(z_k-z^\star)+o(\|z_k-z^\star\|) = M_{\mathrm{PD}}^{k+1-K}(z_K-z^\star)+\sum\nolimits_{j=K}^{k}M_{\mathrm{PD}}^{k-j}o(\|z_j-z^\star\|). \quad \text{(A.19)}$$

Since $z_k \to z^\star$, then take $k$ to the limit we get

$$\lim_{k\to\infty}\sum\nolimits_{j=K}^{k}M_{\mathrm{PD}}^{\infty-j}o(\|z_j-z^\star\|) = \lim_{k\to\infty}\big((z_{k+1}-z^\star)-M_{\mathrm{PD}}^{k+1-K}(z_K-z^\star)\big) = -M_{\mathrm{PD}}^\infty(z_K-z^\star), \quad \text{(A.20)}$$

which holds for any $k\ge K$, *i.e.* $0 = M_{\mathrm{PD}}^\infty(z_k-z^\star)+\lim_{l\to\infty}\sum_{j=k}^{l}M_{\mathrm{PD}}^{l-j}o(\|z_j-z^\star\|)$. Now back to (3.10), for any $k\ge K$,

$$\begin{aligned}
&z_{k+1}-z^\star \\
&= M_{\mathrm{PD}}(z_k-z^\star)+o(\|z_k-z^\star\|) \\
&= M_{\mathrm{PD}}(z_k-z^\star)-M_{\mathrm{PD}}^\infty(z_k-z^\star)+o(\|z_k-z^\star\|)-\lim_{l\to\infty}\sum\nolimits_{j=k}^{l}M_{\mathrm{PD}}^{l-j}o(\|z_j-z^\star\|) \\
&= (M_{\mathrm{PD}}-M_{\mathrm{PD}}^\infty)(z_k-z^\star)+(\mathbf{Id}-\lim_{l\to\infty}M_{\mathrm{PD}}^{l-k})o(\|z_k-z^\star\|)-\lim_{l\to\infty}\sum\nolimits_{j=k+1}^{l}M_{\mathrm{PD}}^{l-j}o(\|z_j-z^\star\|) \\
&= (M_{\mathrm{PD}}-M_{\mathrm{PD}}^\infty)(z_k-z^\star)+(\mathbf{Id}-M_{\mathrm{PD}}^\infty)o(\|z_k-z^\star\|)+M_{\mathrm{PD}}^\infty(z_{k+1}-z^\star),
\end{aligned} \quad \text{(A.21)}$$

where the following equivalence is applied, given finite $k\ge K$

$$M_{\mathrm{PD}}-\lim_{l\to\infty}M_{\mathrm{PD}}^{l-k} = M_{\mathrm{PD}}-M_{\mathrm{PD}}^\infty.$$

Then for $M_{\mathrm{PD}}-M_{\mathrm{PD}}^\infty$, we have

$$M_{\mathrm{PD}}-M_{\mathrm{PD}}^\infty = (M_{\mathrm{PD}}-M_{\mathrm{PD}}^\infty)(\mathbf{Id}-M_{\mathrm{PD}}^\infty).$$

Therefore, for (A.21), move $M_{\text{PD}}^{\infty}(\boldsymbol{z}_{k+1} - \boldsymbol{z}^\star)$ to the other side, we get

$$(\mathbf{Id} - M_{\text{PD}}^{\infty})(\boldsymbol{z}_{k+1} - \boldsymbol{z}^\star) = (M_{\text{PD}} - M_{\text{PD}}^{\infty})(\mathbf{Id} - M_{\text{PD}}^{\infty})(\boldsymbol{z}_k - \boldsymbol{z}^\star) + o((\mathbf{Id} - M_{\text{PD}}^{\infty})\|\boldsymbol{z}_k - \boldsymbol{z}^\star\|), \tag{A.22}$$

which proves (3.15).

(ii) When $R, J^*$ are locally polyhedral around the solution pair $(x^\star, v^\star)$, then the small $o$-term vanishes, and (A.20) implies that for any $k \geq K$

$$M_{\text{PD}}^{\infty}(\boldsymbol{z}_k - \boldsymbol{z}^\star) = 0,$$

*i.e.*, $\boldsymbol{z}_k - \boldsymbol{z}^\star$ belongs to the kernel of $M_{\text{PD}}^{\infty}$. (3.16) then follows. □

**Proof of Theorem 3.11.**

(i) Define $d_k \stackrel{\text{def}}{=} (\mathbf{Id} - M_{\text{PD}}^{\infty})(\boldsymbol{z}_k - \boldsymbol{z}^\star)$ and $\psi_k = o(d_k)$, then from (A.22), for $k \geq K$

$$d_{k+1} = (M_{\text{PD}} - M_{\text{PD}}^{\infty})^{k+1-K} d_K + \sum\nolimits_{j=K}^{k} (M_{\text{PD}} - M_{\text{PD}}^{\infty})^{k-j} \psi_j \tag{A.23}$$

Since the spectral radius $\rho(M_{\text{PD}} - M_{\text{PD}}^{\infty}) < 1$, then from the spectral radius formula, given any $\rho \in ]\rho(M_{\text{PD}} - M_{\text{PD}}^{\infty}), 1[$, there exists a constant $C$ such that, for any $k \in \mathbb{N}$

$$\|(M_{\text{PD}} - M_{\text{PD}}^{\infty})^k\| \leq C\rho^k.$$

Therefore, from (A.23), we get

$$\begin{aligned} \|d_{k+1}\| &\leq \|(M_{\text{PD}} - M_{\text{PD}}^{\infty})^{k+1-K} d_K + \sum\nolimits_{j=K}^{k} (M_{\text{PD}} - M_{\text{PD}}^{\infty})^{k-j} \psi_j\| \\ &\leq \|M_{\text{PD}} - M_{\text{PD}}^{\infty}\|^{k+1-K} \|d_K\| + \sum\nolimits_{j=K}^{k} \|M_{\text{PD}} - M_{\text{PD}}^{\infty}\|^{k-j} \|\psi_j\| \\ &\leq C\rho^{k+1-K}\|d_K\| + C \sum\nolimits_{j=K}^{k} \rho^{k-j} \|\psi_j\|. \end{aligned} \tag{A.24}$$

Together with the fact that $\psi_j = o(\|d_j\|)$ leads to the claimed result.

(ii) When $J$ and $R^*$ are locally polyhedral, then $o(\|\boldsymbol{z}_k - \boldsymbol{z}^\star\|)$ vanishes, and from (3.16) we have

$$\boldsymbol{z}_{k+1} - \boldsymbol{z}^\star = (M_{\text{PD}} - M_{\text{PD}}^{\infty})(\boldsymbol{z}_k - \boldsymbol{z}^\star) = (M_{\text{PD}} - M_{\text{PD}}^{\infty})^{k-K}(\boldsymbol{z}_K - \boldsymbol{z}^\star).$$

Then applying the spectral formula leads to the result (3.18). If $M_{\text{PD}}$ is normal, then it converges linearly to $M_{\text{PD}}^{\infty}$ at the optimal rate $\rho = \rho(M_{\text{PD}} - M_{\text{PD}}^{\infty}) = \|M_{\text{PD}} - M_{\text{PD}}^{\infty}\|$. Combining all this then entails

$$\begin{aligned} \|\boldsymbol{z}_{k+1} - \boldsymbol{z}^\star\| &\leq \|(M_{\text{PD}} - M_{\text{PD}}^{\infty})^{k-K}\|(\boldsymbol{z}_K - \boldsymbol{z}^\star) \\ &= \|M_{\text{PD}} - M_{\text{PD}}^{\infty}\|^{k-K}\|\boldsymbol{z}_K - \boldsymbol{z}^\star\| \\ &= \rho^{k-K}\|\boldsymbol{z}_K - \boldsymbol{z}^\star\|, \end{aligned}$$

and we conclude the proof. □

**Proof of Corollary 5.1.** Owing to the result of [19], condition (5.7) guarantees the convergence of the algorithm.

(i) the identification result follows naturally from Theorem 3.3.

(ii) The result follows Proposition 3.6, Corollary 3.10 and Theorem 3.11. First, for the update of $x_k$ of (5.1), we have

$$x_{k+1} - x^\star = W_{\overline{R}}\overline{H}_F(x_k - x^\star) - \gamma_R W_{\overline{R}} \sum\nolimits_i \overline{L}_i^*(v_{i,k} - v_i^\star) + o(\|x_k - x^\star\| + \gamma_R \sum\nolimits_i \|L_i\|\|v_{i,k} - v_i^\star\|). \tag{A.25}$$

Then the update of $v_{i,k+1}$, for each $i = 1, ..., m$, similar to (A.14), we get

$$\begin{aligned} v_{i,k+1} - v_i^\star &= \big(W_{\overline{J_i^*}}\overline{H}_{G_i^*} - (1+\theta)\gamma_{J_i}\gamma_R W_{\overline{J_i^*}}\overline{L}_i W_{\overline{R}}\overline{L}_i^*\big)(v_{i,k} - v_i^\star) \\ &\quad + \big(2\gamma_{J_i} W_{\overline{J_i^*}}\overline{L}_i W_{\overline{R}}\overline{H}_F - \gamma_{J_i} W_{\overline{J_i^*}}\overline{L}_i\big)(x_k - x^\star) \\ &\quad + o(\|x_k - x^\star\| + \gamma_R \sum\nolimits_i \|L_i\|\|v_{i,k} - v_i^\star\|) + o(\|v_{i,k} - v_i^\star\| + \gamma_{J_i}\|L_i\|\|x_k - x^\star\|). \end{aligned} \tag{A.26}$$

Now consider the small $o$-terms. For the 2 small $o$-terms in (A.10) and (A.13). First, let $a_0, a_1, ..., a_m$ be $m+1$ constants, then we have

$$\sum_{i=0}^{m}|a_i| = \sqrt{(\sum_{i=0}^{m}|a_i|)^2} \leq \sqrt{(m+1)\sum_{i=0}^{m}|a_i|^2} = \sqrt{m+1}\|(a_0, ..., a_m)^T\|.$$

Denote $b = \max\{1, \sum_i \sigma_i\|L_i\|, \gamma_R\|L_1\|, ..., \gamma_R\|L_m\|\}$, then

$$\begin{aligned}&\sum_i \left(\|v_{i,k} - v_i^\star\| + \sigma_i\|L_i\|\|x_k - x^\star\|\right) + \left(\|x_k - x^\star\| + \gamma_R\sum_i \|L_i\|\|v_{i,k} - v_i^\star\|\right) \\ &\leq 2b(\|x_k - x^\star\| + \sum_i\|v_{i,k} - v_i^\star\|) \leq 2b\sqrt{m+1}\|\boldsymbol{z}_k - \boldsymbol{z}^\star\|.\end{aligned}$$

Combining this with (A.25) and (A.26), and ignoring the constants of the small $o$-term, we have that the fixed-point iteration (5.3) is equivalent to

$$\boldsymbol{z}_{k+1} - \boldsymbol{z}^\star = M_{\mathrm{PD}}(\boldsymbol{z}_k - \boldsymbol{z}^\star) + o(\|\boldsymbol{z}_k - \boldsymbol{z}^\star\|).$$

The rest of the proof follows the proof of Theorem 3.11. □